\documentclass[10pt]{article}
\usepackage{amsmath,amssymb,amsthm,amsfonts,enumerate,times,epsfig,color,wrapfig,graphicx,tikz,subfig,pdfpages,mathrsfs,hyperref,cleveref}
\usepackage[numbers]{natbib}
\usepackage[letterpaper,twoside,outer=1.3in,vmargin=1.3in,]{geometry}
\usepackage{hyperref,soul}
\hypersetup{
    colorlinks,
    citecolor=black,
    filecolor=black,
    linkcolor=black,
    urlcolor=black
}
\usepackage{pgfplots}
\pgfplotsset{compat=1.17}

\newtheorem{theorem}{Theorem}[section]
\newtheorem{CO}[theorem]{Corollary}
\newtheorem{LE}[theorem]{Lemma}

\newtheorem{CN}[theorem]{Conjecture}

\newtheorem{EG}[theorem]{Example}
\newtheorem{DE}[theorem]{Definition}
\newcounter{claim_nb}[theorem]
\newtheorem*{claim*}{Claim}

\newcommand{\tr}{\triangle}

\newcommand{\ignore}[1]{}
\newcommand{\conv}{\mathrm{conv}}

\newcommand{\supp}{\mathrm{support}}
\newcommand{\cuboid}{\mathrm{cuboid}}
\newcommand{\core}{\mathrm{core}}

\newcommand{\scr}{\mathrm{SCR}}

\newcommand{\setcore}[1]{\mathrm{setcore}\!\left(#1\right)}

\newcommand{\1}{{\bf 1}}
\newcommand{\0}{{\bf 0}}

\newenvironment{cproof}
{\begin{proof}
[Proof of Claim.]
\vspace{-1.2\parsep}}
{ \end{proof}}

\newif\ifnotes\notesfalse 

\ifnotes
\usepackage{color}
\newcommand{\notename}[2]{{\textcolor{red}{\footnotesize{\bf (#1:} {#2}{\bf ) }}}}

\newcommand{\anote}[1]{{\notename{Ahmad}{#1}}}
\newcommand{\gnote}[1]{{\notename{Gerard}{#1}}}
\newcommand{\onote}[1]{{\notename{Olha}{#1}}}
\newcommand{\snote}[1]{{\notename{Siyue}{#1}}}
\renewcommand{\b}[1]{{\color{blue} #1}}
\renewcommand{\r}[1]{{\color{red} #1}}
\newcommand{\rt}[1]{{\color{red}{\st{#1}}}}
\newcommand{\rb}[2]{\rt{#1}\b{#2}}

\else

\newcommand{\notename}[2]{{}}

\newcommand{\gnote}[1]{}
\newcommand{\anote}[1]{}
\newcommand{\onote}[1]{}
\newcommand{\snote}[1]{}
\renewcommand{\b}[1]{#1}
\renewcommand{\r}[1]{}
\newcommand{\rt}[1]{}
\newcommand{\rb}[2]{#2}

\fi

\title{Lower bounds for cube-ideal set-systems}

\author{Ahmad Abdi\thanks{\b{LSE, Houghton St, London WC2A 2AE, UK, \texttt{a.abdi1@lse.ac.uk} (corresponding author)}} \and G\'{e}rard Cornu\'{e}jols \thanks{\b{Tepper School of Business, CMU, 5000 Forbes Avenue, Pittsburgh, PA 15213, USA}} \and Daniel Dadush \thanks{\b{CWI, Science Park 123, 1098 XG Amsterdam, The Netherlands}} \and Mahsa Dalirrooyfard \thanks{\b{LSE, Houghton St, London WC2A 2AE, UK}}}

\begin{document}
	
\maketitle
	
\begin{abstract}
A set-system $S\subseteq \{0,1\}^n$ is cube-ideal if its convex hull can be described by capacity and generalized set covering inequalities. In this paper, we use combinatorics, convex geometry, and polyhedral theory to give exponential lower bounds on the size of cube-ideal set-systems, and linear lower bounds on their VC dimension. We then provide applications to graph theory and combinatorial optimization, specifically to strong orientations, perfect matchings, dijoins, and ideal clutters, \b{including the Lov\'{a}sz--Plummer conjecture}. 


\noindent \b{\textbf{MSC 2020 Codes.}
03B05, 
05A20, 
05C, 
05C20, 
05C65, 
52A40, 
52B, 
90C27. 
}
\end{abstract}

\section{Introduction}

A set-system $S\subseteq \{0,1\}^n$ is \emph{cube-ideal} if its convex hull can be described by capacity inequalities ${\bf 0}\leq x\leq \1$ and \emph{generalized set covering (GSC)} inequalities, which are of the form $\sum_{i\in I}x_i+\sum_{j\in J}(1-x_j)\geq 1$ for disjoint subsets $I,J\subseteq [n]$. 

Cube-ideal set-systems form a rich class with examples coming from different corners of mathematics\rb{.}{ such as integer programming and graph theory.} For example, if the Hamming graph of $\{0,1\}^n\setminus S$ has degree at most $2$, then $S$ is cube-ideal~\cite{Cornuejols16}. A second example is the cycle space of any graph~\cite{Barahona86}, while strongly connected re-orientations of any digraph give a third class of examples~\cite{Edmonds77}. More recent examples have also been found, see \cite{Abdi-res} for instance. 

In propositional logic, a cube-ideal set-system corresponds to the solutions of a Boolean formula in clausal normal form whose linearization intersected with the unit hypercube forms an integral polytope~\cite{Hooker88,Hooker96}. In {combinatorial optimization}, cube-idealness plays a key role in the study of a basic class of objects known as \emph{ideal clutters}~\cite{Cornuejols94,Cornuejols01}. More specifically, several important conjectures about ideal clutters can be reduced to or equivalently reformulated in terms of cube-ideal set-systems~\cite{Abdi-cuboids,Abdi-kwise}.

In this paper, we give exponential lower bounds on the size of cube-ideal set-systems, and linear lower bounds on their VC dimension. We then provide applications to strong orientations, perfect matchings, dijoins, and ideal clutters. To elaborate, we need two definitions. 

Let $S\subseteq \{0,1\}^n$ be a set-system. The \emph{connectivity of $S$}, denoted by $\lambda(S)$, is the minimum number of variables used in a GSC inequality valid for \b{the convex hull of $S$, denoted by} $\conv(S)$; if $S=\{0,1\}^n$ then $\lambda(S):=+\infty$.\footnote{This notion is closely related to the recently defined notion of \emph{notch} -- more precisely, when connectivity is finite, it is equal to the notch minus $1$~\cite{Benchetrit18}.} \b{If $\lambda(S)=1$, then by coordinate-fixing we can move to a lower-dimensional cube-ideal set-system of the same size; by doing this repeatedly, we may assume that $\lambda(S)\geq 2$.} \r{[line break]}

Denote by $H:[0,\frac12]\to [0,1]$ the \emph{binary entropy} function, defined as $H(\varepsilon) := -\varepsilon \log_2(\varepsilon) - (1 - \varepsilon) \log_2(1-\varepsilon)$ for $\varepsilon > 0$, and $H(0):=0$. Note that $H$ is a strictly increasing continuous function, $H(1/3) \approx 0.9183$, and $H(1/2)=1$. 
\b{The binary entropy function is useful for estimating the number of subsets of a given size. For example, ${n\choose \alpha n}$ is between $2^{H(\alpha)n}/(n+1)$ and $2^{H(\alpha)n}$~(\cite{Mitzenmacher17}, Lemma 10.2), though a more refined asymptotic estimate can be found by using Stirling's approximation, where the dependence on $n$ in the denominator is proportional to $\sqrt{n}$~(see \cite{Galvin14}, \S3.1).} \r{[line break]}

We prove the following lower bound.

\begin{theorem}[proved in \S\ref{sec:large-subcube}]\label{cube-ideal-size-c3}
Let $S\subseteq \{0,1\}^n$ be a cube-ideal set-system with connectivity $\lambda$, where $\lambda\geq 3$. Then $|S|\geq 2^{(1-H(1/\lambda))n}$. 
\end{theorem}

Observe that $1-H(1/\lambda)\geq 1-H(1/3)\geq 0.0817$, so the theorem gives a lower bound of $2^{0.0817n}$ on the size of every cube-ideal set-system $S\subseteq \{0,1\}^n$ with connectivity at least $3$. {We shall see that not only does $\conv(S)$ have exponentially many $0,1$ points, but so do certain high-dimensional faces of it}. 
\rt{To this end, a \emph{set covering (SC)} inequality is an inequality of the form $\sum_{i\in I} x_i\geq 1$ for some subset $I\subseteq [n]$. For an integer $k\geq 2$, a \emph{$k$-SC inequality} is an SC inequality with exactly $k$ variables.}

\begin{theorem}[proved in \S\ref{sec:ideal-c3-core}]\label{cube-ideal-faces-lower-bound}
For every integer $\lambda\geq 3$ and constant $\beta\in (0,1]$, there is a constant $\theta:=\b{\theta_{\ref{cube-ideal-faces-lower-bound}}}(\lambda,\beta)>0$ such that the following statement holds: 

Let $S\subseteq \{0,1\}^n$ be a cube-ideal set-system with connectivity $\lambda${, and let $F$ be the minimal face of $\conv(S)$ containing some $x\in \mathbb{R}^n$, where $x_i\in \left\{\frac{1}{\lambda},1-\frac{1}{\lambda}\right\}$ for each $i\in [n]$. If $F$ has dimension at least $\beta n$, then $|S\cap F|\geq e^{\theta n}$.}
\end{theorem}

What about cube-ideal set-systems $S\subseteq \{0,1\}^n$ with connectivity $2$? 
Note that neither \Cref{cube-ideal-size-c3} nor \Cref{cube-ideal-faces-lower-bound} include the $\lambda=2$ case. One explanation for this disparity is that, in contrast to the $\lambda\geq 3$ case, a cube-ideal set-system with connectivity $2$ can have size linear in $n$ even if $\conv(S)$ is full-dimensional, as we shall see in \Cref{eg:connectivity-2-linear} in \S\ref{sec:ideal}. Nonetheless, for cube-ideal set-systems of connectivity $2$, we can still provide a lower bound on the size of $S$ that is exponential not in $n$ but in the dimension of the minimal face of $\conv(S)$ containing $\frac12 \1$. To elaborate, for an integer $k\geq 2$, a \emph{$k$-GSC inequality} is a GSC inequality that involves exactly $k$ variables.

\begin{DE}\label{DE:cover-graph}
Let $S\subseteq \{0,1\}^n$ be a set-system with connectivity at least $2$. The \emph{$2$-cover graph of $S$}, denoted $G(S)$, is the graph on vertex set $[n]$ with an edge $\{i,j\}$ for every pair $i,j$ of indices such that \b{at least} one of the $2$-GSC inequalities $$\begin{array}{lclclc}
&x_i &+& (1-x_j)&\geq& 1\\ 
&(1-x_i) &+& x_j &\geq& 1\\ 
&x_i &+& x_j&\geq& 1\\
&(1-x_i)&+&(1-x_j)&\geq& 1
\end{array}$$
 is valid for $\conv(S)$. 
The \emph{core of $S$}, denoted $\core(S)$, is the set of points in $S$ that satisfy every $2$-GSC inequality valid for $\conv(S)$ at equality. A GSC inequality $\sum_{i\in I}x_i + \sum_{j\in J}(1-x_j)\geq 1$ valid for $\conv(S)$ is \emph{rainbow} if $I\cup J$ intersects every connected component of $G(S)$ at most once. 
\end{DE}

Observe that every rainbow inequality involves at least $3$ variables. 

When $S$ is a cube-ideal set-system, we shall prove in \S\ref{sec:2-cover-graph} that $G(S)$ is the comparability graph of a preorder, and $\core(S)$ is also a cube-ideal set-system. We shall also prove the following lower bound on $|S|$.

\begin{theorem}[proved in \S\ref{sec:2-cover-graph}]\label{cube-ideal-size}
Let $S\subseteq \{0,1\}^n$ be a cube-ideal set-system with connectivity at least $2$, and fix an inequality description of $\conv(S)$ comprised of capacity and GSC inequalities. Let $d$ be the number of connected components of $G(S)$, and let $\kappa$ be the minimum number of variables used in a rainbow inequality in the fixed description for $\conv(S)$. Then $|\core(S)|\geq 2^{(1-H(1/\kappa))d}$.
\end{theorem}

Note that $\kappa:=+\infty$ if there is no rainbow inequality in the fixed description. Note further that if $S$ has connectivity at least $3$, then every valid GSC inequality is rainbow, so $\kappa$ is simply the connectivity $\lambda$ of $S$\rb{, so}{. Consequently,} \Cref{cube-ideal-size} extends \Cref{cube-ideal-size-c3}.

{In \Cref{cube-ideal-size}, $d$ is the dimension of the minimal face $F$ of $\conv(S)$ containing $\frac12 \1$, and $\core(S) = S\cap F$. The latter follows from the fact that $F$ is the face of $\conv(S)$ obtained by setting all $2$-GSC inequalities at equality. To see the former, note that for $x\in F$, fixing $x_i$ determines the value of $x_j$ for all other vertices $j$ in the same connected component as $i$, so there is at most one degree of freedom within each connected component of $G(S)$. On the other hand, there is at least one degree of freedom within each connected component of $G(S)$, since otherwise $F$ satisfies $x_i=f$ for some coordinate $i\in [n]$ and $f\in [0,1]$, which is a contradiction as $S$ is cube-ideal and $\frac{1}{2} \1\in F$. Consequently, $F$ has dimension $d$.}

In \rb{this light}{light of the above}, we conjecture the following \rb{generalization}{extension} of \Cref{cube-ideal-faces-lower-bound}, in which the condition lower bounding the dimension of the face has been dropped, though the lower bound on the number of $0,1$ points in the face depends instead on this dimension.

\begin{CN}\label{cube-ideal-face-CN}
For every integer $\lambda\geq 3$, there is a constant $\theta \b{:=\theta_{\ref{cube-ideal-face-CN}}(\lambda)}>0$ such that the following statement holds: 

Let $S\subseteq \{0,1\}^n$ be a cube-ideal set-system with connectivity $\lambda$\rt{ such that every variable appears in a valid $\lambda$-SC inequality}. Let $d$ be the dimension of the minimal face $F$ of $\conv(S)$ containing $\frac{1}{\lambda} \1$. Then $|S\cap F|\geq e^{\theta d}$.
\end{CN}

\subsection{Lower bounding the VC dimension of cube-ideal set-systems}

\Cref{cube-ideal-size-c3} may be used to lower bound the `VC dimension' of cube-ideal set-systems with connectivity at least $3$. To this end, let $S\subseteq \{0,1\}^n$ be a set-system, and let $I\subseteq [n]$. The \emph{projection of $S$ onto $I$} is the set-system $\{p_I:p\in S\}$, where $p_I\in \{0,1\}^I$ denotes the restriction of $p$ to the index set $I$. The \emph{Vapnik-Chervonenkis (VC) dimension} of $S$ is the largest integer $d$ such that $S$ has $\{0,1\}^d$ as a projection \b{onto some index set $I$}; we set $d:=0$ if $S=\emptyset$~\cite{Vapnik71}. We shall use \Cref{cube-ideal-size-c3} to prove the following lower bound, where $H^{-1}:[0,1]\to [0,\frac12]$ denotes the inverse of the binary entropy function.

\begin{theorem}[proved in \S\ref{sec:large-subcube}]\label{VC-dim-1}
Let $S\subseteq \{0,1\}^n$ be a cube-ideal set-system with connectivity $\lambda$, where $\lambda\geq 3$. Then $S$ has VC dimension at least $f(\lambda)\cdot n$ where $f(\lambda) = H^{-1}(1-H(1/\lambda))$.
\end{theorem}

The function $f:(2,\infty)\to [0,1]$ is increasing (see \Cref{fig:f-plot}) and $f(3)\geq 0.01013$. In particular, every cube-ideal set-system $S\subseteq \{0,1\}^n$ with connectivity at least $3$ has VC dimension at least $0.01013n$. \b{Based on observations in two special cases of cube-ideal set-systems,} \rb{We}{we} conjecture that \Cref{VC-dim-1} can be strengthened as follows. 

\begin{CN}\label{VC-dim-3}
Let $S\subseteq \{0,1\}^n$ be a cube-ideal set-system with connectivity $\lambda$, where $\lambda\geq 3$. Then $S$ has VC dimension at least $h(\lambda)\cdot n + 1$ where $h(\lambda) = 1-2/\lambda$. 
\end{CN}

An intermediate result would be to prove the above for the function $g:(2,\infty)\to [0,1]$ defined as $g(\lambda) := 1-H(1/\lambda)$. Given that every set-system of VC dimension $d$ has size at least $2^d$, this result would imply \Cref{cube-ideal-size-c3}. 

\Cref{VC-dim-3} holds for an important class of examples from graphs. Let $\lambda\geq 3$ be an integer, let $G=(V,E)$ be a $\lambda$-edge-connected graph, and let $S\subseteq \{0,1\}^E$ be the \emph{cycle space} of $G$, that is, $S$ consists of the incidence vector of every edge-subset of $G$ where every vertex has even degree. It is known that $S$ is cube-ideal with connectivity $\lambda$~[\cite{Barahona86}, \b{(3.5)}]. Observe that $S$ is a $GF(2)$-vector space of $GF(2)$-rank $|E|-|V|+1$. This implies that $S$ has VC dimension $|E|-|V|+1$, and $|S|=2^{|E|-|V|+1}$. As $G$ has minimum degree at least $\lambda$, it follows that $|E|-|V|+1\geq (1-2/\lambda)|E|+1$, where the inequality holds at equality if, and only if, $G$ is $\lambda$-regular. Thus, \Cref{VC-dim-3} holds for this example, and the lower bound can be tight. In the next subsection, we see that this conjecture also holds for another class of examples from graphs.

In \S\ref{sec:ideal}, we shall use the famous \emph{width-length} inequality for ideal clutters to prove \Cref{VC-dim-3} for up-monotone set-systems, in fact, we will show that they have VC dimension at least $(1-1/\lambda)n$. 

\begin{figure}
\centering
\begin{tikzpicture}
\begin{axis}[
    width=12cm, height=9cm,
    xmin=2, xmax=20,
    ymin=0, ymax=0.9,
    xlabel={$\lambda$},
    xtick={2,4,6,8,10,12,14,16,18,20},
    ytick={0,0.1,0.2,0.3,0.4,0.5,0.6,0.7,0.8,0.9},
    grid=none,
    axis lines=box,
    every axis x label/.style={at={(ticklabel* cs:0.5)},anchor=north},
    thick,
    legend style={
      at={(1.02,0.5)},
      anchor=west,
      font=\small,
      draw=none,
      fill=none,
      cells={anchor=west}
    },
]

\addplot[black, solid, very thick, smooth]
  table[x=x, y=f] {figures/fdata.dat};
\addlegendentry{$f(\lambda)=H^{-1}\!\bigl(1-H(1/\lambda)\bigr)$}

\addplot[black, dashed, ultra thick, smooth, samples=100, domain=3:20]
  {1 - ( -(1/x)*ln(1/x)/ln(2) - (1-1/x)*ln(1-1/x)/ln(2) )};
\addlegendentry{$g(\lambda)=1-H(1/\lambda)$}

\addplot[black, dotted, ultra thick, smooth, samples=100, domain=3:20]
  {1 - 2/x};
\addlegendentry{$h(\lambda)=1-2/\lambda$}

\end{axis}
\end{tikzpicture}
\caption{\b{Graphs of $f,g,h$ (from bottom to top) over the domain $[3,20]$.}}
\label{fig:f-plot}
\end{figure}

\subsection{Applications to graph theory and combinatorial optimization}

\r{Reviewer comment: I'm not entirely sure how relevant the applications of the developed results really are, but I still like the general approach and see the main value in the concept itself.}

Our results\b{, in particular \Cref{cube-ideal-size-c3} and \Cref{cube-ideal-faces-lower-bound},} have a number of applications to \b{counting} objects that are of importance in graph theory and combinatorial optimization, including perfect matchings, strong orientations, dijoins, and ideal clutters. We shall briefly explain two of these here, and postpone the remaining two to \S\ref{sec:ideal} and \S\ref{sec:apps}. \b{More specifically, in \S\ref{sec:apps}, we prove part of the Lov\'{a}sz--Plummer conjecture for a basic class of regular `matching-covered' graphs.}\medskip

Let $G=(V,A,E)$ be a mixed graph with arcs $A$, which are directed, and edges $E$, which are not yet directed. A \emph{$G$-strong orientation} is an orientation $\overrightarrow{E}$ of the edges such that the digraph $(V,A\cup \overrightarrow{E})$ is strongly connected. A blocking notion is that of a \emph{pseudo dicut}, which is a cut of the form $\delta_G(U)\subseteq A\cup E$, where $U\subset V$, $U\neq \emptyset$, and $\delta_A^-(U)=\emptyset$. Observe that an orientation $\overrightarrow{E}$ is $G$-strong if, and only if, in every pseudo dicut $\delta_G(U)$ at least one edge of $E$ is oriented to enter $U$.

\rt{Suppose every pseudo dicut of $G=(V,A,E)$ contains at least two edges from~$E$.} Let $\overrightarrow{E}$ be an arbitrary reference orientation of $E$, and denote by $\scr(G;\overrightarrow{E})$ the set of vectors $x\in \{0,1\}^{\overrightarrow{E}}$ such that the re-orientation of $\overrightarrow{E}$ obtained after flipping the arcs in $\{a\in \overrightarrow{E}:x_a=1\}$ is $G$-strong. We shall see that $\scr(G;\overrightarrow{E})$ is a cube-ideal set-system, and thus use our results to obtain the following consequence.

\begin{theorem}[proved in \S\ref{sec:apps}]\label{strong-orientations-exp}
Let $G=(V,A,E)$ be a mixed graph where every pseudo dicut contains at least $\lambda$ edges from~$E$, where $\lambda\geq 3$. Then the number of $G$-strong orientations is at least $2^{(1-H(1/\lambda))|E|}$.
\end{theorem}

This is the first exponential lower bound on the number of $G$-strong orientations of such a mixed graph, as far as we know, \b{and sheds some light on a poorly understood parameter}. \r{[line break]}

\b{In contrast, when $A=\emptyset$, counting the number of $G$-strong orientations has been extensively studied. While computing this number is $\#P$-hard~\cite{Jaeger90}, it is precisely the valuation of the Tutte polynomial at $(0,2)$~(see \cite{EM22}), and plays a key role in the well-known Merino--Welsh conjecture~\cite{Merino99}.}
\rb{When $A=\emptyset$}{In this case}, the lower bound in \Cref{strong-orientations-exp} is easy to obtain and can in fact be improved to $2^{(1-\frac{2}{\lambda})|E|+1}$ by using the theory of ear decompositions as follows. \b{Given a graph $H$, an \emph{$H$-path} is a path that is edge-disjoint from $H$, and shares only its end-vertices with $H$.} Every $2$-vertex-connected block of $G$ can be constructed from a cycle by successively adding $H$-paths to graphs $H$ already constructed; the cycle and the paths are referred to as \emph{ears}. It can be readily checked that the number of ears is exactly {$|E|-|V|+1$, as $G$ is connected}~[\b{\cite{Whitney32}, Theorem 18, see also} \cite{Diestel25}, \rb{Chapter 3}{Proposition 3.1.3}]. Note that by orienting each ear so that it becomes a directed cycle or path, one obtains a strong orientation of $G$. The number of such orientations is $2^{|E|-|V|+\b{1}}$. As $G$ is $\lambda$-edge-connected, every vertex has degree at least $\lambda$, so $|E|-|V|+\b{1}\geq (1-\frac{2}{\lambda})|E|+1$, implying in turn that the number of strong orientations of $G$ is at least $2^{(1-\frac{2}{\lambda})|E|+1}$. In fact, we just proved a lower bound of $|E|-|V|+\b{1}\geq (1-\frac{2}{\lambda})|E|+1$ on the VC dimension of $\scr(G;\overrightarrow{E})$, thus confirming \Cref{VC-dim-3} for $\scr(G;\overrightarrow{E})$ when $A=\emptyset$.
\medskip

Let us exhibit another application. A \emph{bipartite digraph} is a digraph $D=(V,A)$ where every vertex is a source or a sink. A \emph{dicut} is a cut of the form $\delta^+(U)\subseteq A$ where $\delta^-(U)=\emptyset$, $U\subset V$, $U\neq \emptyset$. A blocking notion is that of a \emph{dijoin}, which is a subset of $A$ that intersects every dicut at least once. An important open problem in combinatorial optimization is Woodall's conjecture, which states that the minimum size of a dicut, say $\tau$, is equal to the maximum number of arc-disjoint dijoins~\cite{Woodall78}. It suffices to prove this conjecture for bipartite digraphs where every sink has degree~$\tau$~\cite{Abdi-dijoins}. \rb{We shall prove the following theorem as another consequence of our results.}{Observe that in such a digraph, and in any assumed collection of $\tau$ arc-disjoint dijoins, every dijoin is (inclusion-wise) minimal and intersects every minimum dicut exactly once. Motivated by Woodall's conjecture and this observation, we shall prove an exponential lower bound on the number of such dijoins.}

\begin{theorem}[proved in \S\ref{sec:apps}]\label{tight-dijoins-exponential}
For every integer $\tau\geq 3$, there is a constant $\theta\b{:=\theta_{\ref{cube-ideal-faces-lower-bound}}(\tau,\frac13)}>0$ such that the following statement holds: 

Let $D=(V,A)$ be a bipartite digraph where every sink has degree $\tau$, and every dicut has size at least $\tau$. Then the number of minimal dijoins of $D$ intersecting every minimum dicut exactly once is at least $e^{\theta |A|}$.
\end{theorem}

\section{$3$-connected cube-ideal set-systems}\label{sec:large-subcube}

In this section, we prove \Cref{cube-ideal-size-c3}, which states that if $S\subseteq \{0,1\}^n$ is a cube-ideal set-system with connectivity $\lambda\geq 3$, then $|S|\geq 2^{(1-H(1/\lambda))n}$. We will then use this result to prove a lower bound on the VC dimension of such set-systems, namely \Cref{VC-dim-1}. The following observation is crucial.

\begin{LE}\label{cube-ideal-large-subcube}
Let $S\subseteq \{0,1\}^n$ be a cube-ideal set-system with connectivity $\lambda$, where $\lambda\geq 3$. Then $\conv(S)\supseteq \big[\frac{1}{\lambda},1-\frac{1}{\lambda}\big]^n$.
\end{LE}
\begin{proof}
As $S$ is cube-ideal, $\conv(S)$ is described by capacity and GSC inequalities. As $S$ has connectivity $\lambda$, every GSC inequality in the description must involve at least $\lambda$ variables. Thus, any point in $\big[\frac{1}{\lambda},1-\frac{1}{\lambda}\big]^n$ satisfies the inequality description of $\conv(S)$, implying the desired containment.
\end{proof}

We will show that the containment above \rb{alone implies that $|S|$ must be}{is the culprit for $|S|$ being} exponentially large. We need the following well-known inequality.

\begin{LE}[see \cite{Galvin14}, Theorem 3.1, also \cite{Mitzenmacher17}, \S 10.2]\label{entropy}
For \rb{all integers}{every integer} $n\geq 1$ and \b{any real number} $\lambda\geq 2$, the number of subsets of $[n]$ of size at most $n/\lambda$ is at most $2^{H(1/\lambda)n}$.
\end{LE}

Given $S\subseteq \{0,1\}^n$ and $q\in \{0,1\}^n$, to \emph{twist $S$ by $q$} is to replace $S$ by $S\tr q:=\{p\tr q:p\in S\}$, where the second $\tr$ denotes coordinate-wise addition modulo $2$. Take a coordinate $i\in [n]$. Denote by $e_i$ the $i\textsuperscript{th}$ unit vector of appropriate dimension. To \emph{twist coordinate $i$ of $S$} is to replace $S$ by $S\tr e_i$. 

We are now ready to prove the following.

\begin{theorem}\label{large-subcube}
Take an integer $\lambda\geq 3$. Let $S\subseteq \{0,1\}^n$ be a set-system such that $\conv(S)\supseteq \big[\frac{1}{\lambda},1-\frac{1}{\lambda}\big]^n$. Then $|S|\geq 2^{(1-H(1/\lambda))n}$. 
\end{theorem}
\begin{proof}
For each $w\in \{-1,+1\}^n$, let $x[w]$ be a point in $S$ which maximizes $w^\top x$. 
As $\conv(S)\supseteq \big[\frac{1}{\lambda},1-\frac{1}{\lambda}\big]^n$, we have the inequality below, 
$$w^\top x[w]=\max \left\{w^\top x:x\in \conv(S)\right\}\geq \frac{\lambda-1}{\lambda} |\supp_+(w)|-\frac{1}{\lambda} |\supp_-(w)|,$$ where $\supp_+(w)=\{i\in [n]:w_i=+1\}$ and $\supp_-(w)=\{i\in [n]:w_i=-1\}$. There exists an $x^\star\in S$ such that $$
\left|\left\{w\in \{-1,+1\}^n : w^\top x^\star = w^\top x[w]\right\}\right|\geq \frac{2^n}{|S|}.$$ Let $T:=\left\{w\in \{-1,+1\}^n : w^\top x^\star = w^\top x[w] \right\}$. Note that for any $w\in T$, we have $$
w^\top x^\star\geq \frac{\lambda-1}{\lambda} |\supp_+(w)|-\frac{1}{\lambda} |\supp_-(w)|.
$$
After twisting the coordinates, if necessary, we may assume that $x^\star={\bf 0}$; observe that twisting coordinate $i$ maps $x_i\mapsto 1-x_i$ and $w_i\mapsto -w_i$. This means that for any $w\in T$, we have $$
0\geq \frac{\lambda-1}{\lambda} |\supp_+(w)|-\frac{1}{\lambda} |\supp_-(w)|,
$$ or equivalently, $|\supp_+(w)| \leq n / \lambda$. This implies in turn that $$
T \subseteq \{ w \in \{ -1, 1 \}^n: |\supp_+(w)| \leq n / \lambda \}.
$$ 
Subsequently, $$
|T|\leq \left|\left\{y\in \{0,1\}^n : \1^\top y\leq n/\lambda\right\}\right|\leq 2^{H(1/\lambda)n}
$$ where the rightmost inequality follows from \Cref{entropy}. Putting the inequalities together, we get $$
2^{H(1/\lambda)n}\geq |T|\geq \frac{2^n}{|S|}
$$ implying in turn that $|S|\geq 2^{(1-H(1/\lambda))n}$, as desired.
\end{proof}

We are now ready to prove \Cref{cube-ideal-size-c3}.

\begin{proof}[Proof of \Cref{cube-ideal-size-c3}]
Let $S\subseteq \{0,1\}^n$ be a cube-ideal set-system with connectivity $\lambda$, where $\lambda\geq 3$. By \Cref{cube-ideal-large-subcube}, $\conv(S)\supseteq \big[\frac{1}{\lambda},1-\frac{1}{\lambda}\big]^n$, so $|S|\geq 2^{(1-H(1/\lambda))n}$ by \Cref{large-subcube}, as required.
\end{proof}

Next, we use the above to give a lower bound on the VC dimension of cube-ideal set-systems with connectivity at least $3$. We shall need the following well-known lemma.

\begin{LE}[Sauer-Shelah~\cite{Sauer72,Shelah72}]\label{SS-lemma}
Let $n\geq k\geq 1$ be integers, and let $S\subseteq \{0,1\}^n$ be a set-system of VC dimension at most $k$. Then $|S|$ is at most the number of subsets of $[n]$ of size at most $k$.
\end{LE}

We need the following convenient consequence of the lemma.

\begin{CO}\label{SS-lemma-CO}
Let $\lambda\geq 2$ be an integer, and let $S\subseteq \{0,1\}^n$ be a set-system such that $|S|\geq 2^{H(1/\lambda)n}$. Then $S$ has VC dimension at least $\frac{n}{\lambda}$.
\end{CO}
\begin{proof}
Let $k^\circ:=\frac{n}{\lambda}$ and $k:=\lfloor k^\circ\rfloor$. We have that $|S|\geq 2^{H(1/\lambda)n}\geq 2^{H(k/n)n}$ as $H$ is a strictly increasing function. Furthermore, by \Cref{entropy}, $2^{H(k/n)n}$ is greater than or equal to the number of subsets of $[n]$ of size at most $k$.

If $k^\circ =k$, i.e., $k^\circ$ is an integer, then $|S|$ is strictly larger than the number of subsets of $[n]$ of size at most $k-1$. It therefore follows from (the contrapositive of) \Cref{SS-lemma} that $S$ has VC dimension strictly larger than $k-1$, so $S$ has VC dimension at least $k^\circ$.

Otherwise, $k^\circ>k$, so $|S|\geq 2^{H(1/\lambda)n}>2^{H(k/n)n}$. Subsequently, by  \Cref{entropy}, $|S|$ is strictly larger than the number of subsets of $[n]$ of size at most $k$. It therefore follows from (the contrapositive of) \Cref{SS-lemma} that $S$ has VC dimension strictly larger than $k$, so $S$ has VC dimension at least $k+1>k^\circ$.
\end{proof}

We are ready to prove \Cref{VC-dim-1}.

\begin{proof}[Proof of \Cref{VC-dim-1}]
Let $S\subseteq \{0,1\}^n$ be a cube-ideal set-system with connectivity $\lambda$, where $\lambda\geq 3$. We know from \Cref{cube-ideal-size-c3} that $|S|\geq 2^{(1-H(1/\lambda))n}$. It therefore follows from \Cref{SS-lemma-CO} that $S$ has VC dimension at least $H^{-1}(1-H(1/\lambda)) \cdot n$, as desired.
\end{proof}

\section{Faces of cube-ideal set-systems with connectivity at least $3$}\label{sec:ideal-c3-core}

In this section, we prove \Cref{cube-ideal-faces-lower-bound} \rb{on a certain face}{giving exponential lower bounds for certain high-dimensional faces} of the convex hull of a cube-ideal set-system \b{$S\subseteq \{0,1\}^n$} with connectivity at least $\b{\lambda\geq} 3$. \b{Let us give a summary of the key ideas of the proof. To this end, let $F$ be the minimal face of $\conv(S)$ containing $\frac{1}{\lambda}\1$, say, of dimension $d\geq \beta n$. To give an exponential lower bound on $|S\cap F|$, the proof proceeds by showing that the face $F$ contains a `hypercube-like' polytope that is centrally symmetric around $\frac{1}{\lambda}\1$, and has `large radius' and `few' facets. More specifically,  $F$ will contain an affine slice of a hypercube with constant positive $\ell_\infty$-radius and $O(d)$ facets, namely, for $\varepsilon:=\frac{1}{\lambda}-\frac{1}{\lambda+1}>0$, 
$$F\supseteq \left[\frac{1}{\lambda}-\varepsilon,\frac{1}{\lambda}+\varepsilon\right]^n \cap \mathrm{aff}(F)=:Q,$$ where $\mathrm{aff}(F)$ refers to the affine hull of $F$, and $Q$ is a $d$-dimensional polytope with at most $2n\leq \frac{2}{\beta} d$ facets. At this point, by resorting to a convex geometric tool by Barvinok~\cite{Barvinok13}, we can conclude that $F$ has exponentially many vertices, thus proving \Cref{cube-ideal-faces-lower-bound}.} \rt{Below, $\mathrm{aff}(F)$ refers to the affine hull of $F$.}

\b{We start by proving the inclusion above.}

\begin{LE}\label{LE:cube-ideal-face-large-subcube}
Let $\lambda\geq 3$ be an integer, and let $S\subseteq \{0,1\}^n$ be a cube-ideal set-system with connectivity $\lambda$\rt{ such that if $x_i,i\in [n]$ appears in a valid $\lambda$-GSC inequality, then it appears in a valid $\lambda$-SC inequality}. Let $F$ be the minimal face of $\conv(S)$ containing $\frac{1}{\lambda} \1$. Then {$\big[\frac{1}{\lambda}-\varepsilon,\frac{1}{\lambda}+\varepsilon\big]^n\cap \mathrm{aff}(F) \subseteq F$ for $\varepsilon=\frac{1}{\lambda}-\frac{1}{\lambda+1}$}. 
\end{LE}
\begin{proof}
Let $\mathcal{V}$ be the set of all pairs $(I,J)$ such that $I,J$ are disjoint subsets of $[n]$, and $\sum_{i\in I} x_i+\sum_{j\in J} (1-x_j)\geq 1$ is valid for $\conv(S)$. Such an inequality is tight at $\frac{1}{\lambda} \1$ if, and only if, $J=\emptyset$ and $|I|=\lambda$; here we used the inequality $1-\frac{1}{\lambda}>\frac{1}{\lambda}$ guaranteed by $\lambda\geq 3$. Thus, \begin{align*}
F &= \conv(S)\cap \{ x \b{\in \mathbb{R}^n}: x(I)= 1,\, \forall~(I,\emptyset)\in \mathcal{V} \text{ s.t. } |I|=\lambda \}\\
\mathrm{aff}(F) &= \{ x \b{\in \mathbb{R}^n} : x(I)= 1,\, \forall~(I,\emptyset)\in \mathcal{V} \text{ s.t. } |I|=\lambda \}.
\end{align*}
Let {$p\in \big[\frac{1}{\lambda}-\varepsilon,\frac{1}{\lambda}+\varepsilon\big]^n\cap \mathrm{aff}(F)$}. \b{Note that for each $i\in [n]$, $$\frac{1}{\lambda+1}=\frac{1}{\lambda}-\varepsilon\leq p_i\leq \frac{1}{\lambda}+\varepsilon< \frac{1}{2}<1-\frac{1}{\lambda+1}.$$} 
We need to show that $p\in F$.
To this end, let $(I,J)\in \mathcal{V}$. If $|I|+|J|\geq \lambda+1$, then $p$ satisfies the corresponding GSC inequality as $p\in \big[\frac{1}{\lambda+1},1-\frac{1}{\lambda+1}\big]^n$. If $|I|=\lambda$ and $J=\emptyset$, then $p$ satisfies the corresponding GSC inequality at equality as $p\in \mathrm{aff}(F)$. 
\b{If $|J|\geq 2$, then $\sum_{i\in I} p_i + \sum_{j\in J} (1-p_j)\geq \sum_{j\in J} (1-p_j)>|J|\frac12 \geq 1$.}
Otherwise, {$|I|=\lambda-1$ and $|J|=1$. Then $$\sum_{i\in I} p_i + \sum_{j\in J} (1-p_j)\geq (\lambda-1)\left(\frac{1}{\lambda+1}\right)+1-\frac{1}{\lambda}-\varepsilon = 1+\frac{\lambda}{\lambda+1}-\frac{2}{\lambda}>1,$$ as required.} \rt{Let $k\in J$. By hypothesis, there exists $(K,\emptyset)\in \mathcal{V}$ such that $|K|=\lambda$ and $k\in K$. Given that $p$ satisfies the corresponding GSC inequality at equality, it follows that $1-p_k = \sum_{\ell\in K\setminus k} p_\ell$, where the right-hand side involves $|K|-1=\lambda-1\geq 2$ variables. Subsequently, as $p\in \big[\frac{1}{\lambda+1},1-\frac{1}{\lambda+1}\big]^n$, we have $
\sum_{i\in I} p_i +$ $\sum_{j\in J} (1-p_j) = \sum_{i\in I} p_i + \sum_{\ell\in K\setminus k} p_\ell + \sum_{j\in J\setminus k} (1-p_j) \geq \frac{1}{\lambda+1} \left(|I| + 2 + |J|-1\right) =1,
$ so $p$ satisfies the corresponding GSC inequality, as required.} 
\end{proof}

In what follows, we shall use the inclusion above to lower bound the number of $0,1$ points inside the face $F$. To this end, we need the following theorem.

\begin{theorem}[\cite{Barvinok13}]\label{barvinok}
For every $\alpha,\beta\geq 1$ there is \b{a constant} $\gamma:=\b{\gamma_{\ref{barvinok}}}(\alpha,\beta)>0$ such that the following statement holds: 

Suppose that $P\subset \mathbb{R}^d$ is a polytope containing the set $$\left\{x\in \mathbb{R}^d : -1\leq u_i^\top x\leq 1,\,\forall i\in [m]\right\},$$ where $\|u_i\|_2\leq 1,\,\forall i\in [m]$ and $m\leq \alpha d$. Suppose further that $P$ lies inside the ball $$
\left\{x\in \mathbb{R}^d : \|x\|_2\leq \beta\sqrt{d}\right\}.
$$ Then $P$ has at least $e^{\gamma d}$ vertices.
\end{theorem}

\b{Barvinok~\cite{Barvinok13} gives an explicit formula for the constant in \Cref{barvinok}: Choose any $\varepsilon\in (0,1)$ and a sufficiently large $\rho:=\rho(\alpha,\beta,\varepsilon)>0$ such that the following inequalities hold: 
\begin{align*}
&\alpha\ln\left(1 - \exp\left\{-\frac{\rho^2}{2}\right\}\right) > -\frac{\varepsilon^2}{4}\\
&\gamma_{\ref{barvinok}}(\alpha,\beta):= \frac{(1 - \varepsilon)^2 }{2 \beta^2 \rho^2} \left(1 - \exp\left\{-\frac{\rho^2}{2}\right\}\right) +  \alpha\ln\left(1 - \exp\left\{-\frac{\rho^2}{2}\right\}\right)>0.
\end{align*} 
\Cref{barvinok} is then proved for all sufficiently large $d>d_0(\alpha,\beta,\varepsilon,\rho)$ and for $\gamma:=\gamma_{\ref{barvinok}}(\alpha,\beta)$. Note the slight notational shortcoming, as the definition of $\gamma_{\ref{barvinok}}(\alpha,\beta)$ also depends on the choice of $\varepsilon,\rho$.}

\paragraph{Digression.} We can apply \Cref{barvinok} to get an alternate proof of \Cref{large-subcube} for sufficiently large $n$ but with \rb{an inexplicit}{a weaker} lower bound. More specifically, we let $P = \frac{2\lambda}{\lambda-2}(\conv(S) - \frac12 \1)$. Then $P$ contains $\left[-1,1\right]^n$ and is contained in a Euclidean ball of radius $\frac{\lambda}{\lambda-2}\sqrt{n}$ around the origin. Thus, we can set $\alpha = 1$ and $\beta=\frac{\lambda}{\lambda-2}$\b{, and let $\gamma:=\gamma_{\ref{barvinok}}\big(1,\frac{\lambda}{\lambda-2}\big)$}. \rt{For sufficiently large~$n$, Barvinok~[??] gives an explicit formula for $\gamma:=\gamma(\alpha,\beta)$, which for our choices of $\alpha,\beta$ leads to the following: Choose any $\varepsilon\in (0,1)$ and $\rho>0$ such that the following inequality holds:
$
\ln\left(1 - \exp\left\{-\frac{\rho^2}{2}\right\}\right) > -\frac{\varepsilon^2}{4}.
$
Then for all sufficiently large $n>n_0(\alpha,\beta,\varepsilon,\rho)$ we can choose $\gamma$ as follows:
$
\gamma = \frac{(\lambda-2)^2 (1 - \varepsilon)^2 }{2 \lambda^2 \rho^2} \left(1 - \exp\left\{-\frac{\rho^2}{2}\right\}\right) +  \ln\left(1 - \exp\left\{-\frac{\rho^2}{2}\right\}\right) > 0.
$} 
In contrast, \Cref{large-subcube} gives {$\hat{\gamma}=(1-H(1/\lambda))\ln{2}$} for all \rb{$n$}{$n\geq 1$}. 

\b{At $\lambda=3,10^2,10^6,10^{12}$, the values for $\hat{\gamma}$ are approximately 
$
0.0566330, 0.6371456, 0.6931323,0.6931472$, respectively. However, for $\varepsilon=0.1$ and $\rho=3.462$, the values for $\gamma$ are approximately $
0.0012451, 0.0298718,$ $0.0312065,0.0312067$, respectively. Given $\lambda=10^6$, the values $
\varepsilon\approx 0.089650\in (0,1),\,
\rho\approx 3.524481>0$ maximize $\gamma$ at $\gamma^\star\approx 0.0312814$, according to Matlab's \texttt{fmincon} function within the optimization toolbox; we note that $\hat{\gamma}>20\cdot \gamma^\star$.} 
\medskip

Moving on, we are ready to prove the following lemma, which borrows some ideas from [\cite{Barvinok13}, Corollary~1.3].

\begin{LE}\label{faces-lower-bound}
For every integer $\lambda\geq 2$ and constants $\alpha\in \left(\frac{1}{\lambda+1},\frac{1}{2}\right]$ and $\beta\in (0,1]$, there is a constant \rt{$\theta>0$} $\theta\b{:=\beta\gamma_{\ref{barvinok}}\left(\frac{1}{\beta}, \frac{(1+\alpha)}{\varepsilon\sqrt{\beta}}\right)}>0$ \b{for $\varepsilon=\alpha-\frac{1}{\lambda+1}$,} such that the following statement holds: 

Let $S\subseteq \{0,1\}^n$ be a set-system whose convex hull contains $\alpha \1$ (possibly on its boundary), and let $F$ be the minimal face of the convex hull containing $\alpha \1$. Suppose \rt{that} $F$ has dimension at least $\beta n$, and {it contains
$\big[\frac{1}{\lambda+1},2\alpha-\frac{1}{\lambda+1}\big]^n\cap \mathrm{aff}(F)$
}. 
Then $|S\cap F|\geq e^{\theta n}$.
\end{LE}
\begin{proof}
Let $d:=\dim(F)\geq \beta n$, let $A$ be the affine hull of $F$, and let $R:=S\cap F\subseteq \{0,1\}^n\cap A$. By hypothesis, $\conv(R)\supseteq \big[\frac{1}{\lambda+1},\b{2\alpha}-\frac{1}{\lambda+1}\big]^n\cap A$. Let us shift \rb{$\conv(R),A$}{$\conv(R)$} so that \rb{they both contain}{it contains} the origin. To this end, let $a:=\alpha \1$, $P:=\conv(R)-a$, and $L:=A-a$ which is a linear subspace of dimension $d$ containing $P$.

First we show that $P$ contains a large hypercube-like polytope with $O(d)$ facets. To this end, let $\varepsilon:=\alpha-\frac{1}{\lambda+1}>0$. We claim that $Q:=[-\varepsilon,\varepsilon]^n\cap L\subseteq P$. To see this, observe that $a+y\in \big[\frac{1}{\lambda+1},\b{2\alpha}-\frac{1}{\lambda+1}\big]^n$ for each $y\in [-\varepsilon,\varepsilon]^n$, so $a+y\in \conv(R)$ for all $y\in [-\varepsilon,\varepsilon]^n$ such that $a+y\in A$, implying in turn that $Q\subseteq P$. Observe that $$
P\supseteq Q = \left\{x\in L:-\varepsilon \leq u_i^\top x \leq \varepsilon,\,\forall i\in [n]\right\},
$$ where $u_i\in L$ is the orthogonal projection of the standard unit vector $e_i\in \mathbb{R}^n$ onto $L$, for $i\in [n]$. In particular, $\|u_i\|_2\leq 1$ for all $i\in [n]$. Furthermore, we know that $n\leq \frac{1}{\beta} d$. 

Secondly, we show that $P$ is contained in a Euclidean ball of radius $O(\sqrt{d})$. For every $p\in S$, we have $$
\|p - a\|_2\leq \|p\|_2 + \|a\|_2\leq \sqrt{n} + \alpha \sqrt{n} = (1+\alpha) \sqrt{n} \leq \frac{1+\alpha}{\sqrt{\beta}} \sqrt{d}.
$$ Thus, $$P\subset B:=\left\{x\in L:\|x\|_2 \leq \frac{1+\alpha}{\sqrt{\beta}} \sqrt{d}\right\}.$$

Finally, we apply a linear transformation from $L\subset \mathbb{R}^n$ to $\mathbb{R}^d$. To this end, let $\ell_1,\ldots,\ell_d$ be an orthonormal basis for $L$, and let $M\in \mathbb{R}^{n\times d}$ be the matrix whose columns are $\ell_1,\ldots,\ell_d$. Let $f:L\to \mathbb{R}^d$ be the linear transformation that maps $\ell_i$ to $e_i\in \mathbb{R}^d$ for $i\in [d]$. Note that $f^{-1}(z) = Mz$. 

We now work with $f(P)$ instead of $P$. We have
$$
f(P)\supseteq 
f(Q) = \left\{z\in \mathbb{R}^d : -\varepsilon\leq w_i^\top z \leq \varepsilon, \forall i\in [n]\right\},
$$ where $w_i = M^\top u_i$ for each $i\in [n]$. Observe that $$\|w_i\|_2 \leq \|M^\top\|_2 \|u_i\|_2 = \sqrt{\lambda_{\max}(MM^\top)} \cdot \|u_i\|_2 = \|u_i\|_2\leq 1,$$ where $\|M\|_2$ denotes the spectral norm of $M$, which is equal to the square root of the largest eigenvalue of $MM^\top$; the latter has the same nonzero spectrum as $M^\top M=I_d$, so $\lambda_{\max}(MM^\top)=1$. 

Next we have $$f(P)\subset f(B)\b{=}\left\{z\in \mathbb{R}^d : \|Mz\|_2 \leq \frac{1+\alpha}{\sqrt{\beta}} \sqrt{d}\right\} = \left\{z\in \mathbb{R}^d : \|z\|_2 \leq \frac{1+\alpha}{\sqrt{\beta}} \sqrt{d}\right\},$$ where the second equality follows from $\|Mz\|_2 = \sqrt{z^\top M^\top M z} = \|z\|_2$. 

Subsequently, by applying \Cref{barvinok} to $\varepsilon^{-1} f(P)$, we obtain that $\varepsilon^{-1} f(P)$ has at least $e^{\gamma d}$ vertices, for some constant $\gamma:=\b{\gamma_{\ref{barvinok}}}\left(\frac{1}{\beta}, \frac{\varepsilon^{-1}(1+\alpha)}{\sqrt{\beta}}\right)>0$. 
Observe that $\varepsilon^{-1} f(P)$ has the same number of vertices as $f(P)$ and also~$P$, namely $|R|$. Thus, $|S\cap F|=|R|\geq e^{\gamma d}\geq e^{\gamma \beta n}$, so $\theta:=\gamma \beta>0$ is the desired constant.
\end{proof}

We are now ready to prove \Cref{cube-ideal-faces-lower-bound}.

\begin{proof}[Proof of \Cref{cube-ideal-faces-lower-bound}]
\b{After possibly twisting $S$, we may assume that $x=\frac{1}{\lambda}\1$.} 
It follows from \Cref{LE:cube-ideal-face-large-subcube} that {$\big[\frac{1}{\lambda+1},\frac{2}{\lambda}-\frac{1}{\lambda+1}\big]^n\cap \mathrm{aff}(F) \subseteq F$}. We can therefore apply \Cref{faces-lower-bound} for $\alpha = \frac{1}{\lambda}$ to obtain the result \b{for} $$\b{\theta:=\theta_{\ref{cube-ideal-faces-lower-bound}}(\lambda,\beta):= \beta
\gamma_{\ref{barvinok}}\left(\frac{1}{\beta}, \frac{(\lambda+1)^2}{\sqrt{\beta}}\right)>0.}
$$
\end{proof}

\section{$2$-connected cube-ideal set-systems}\label{sec:2-cover-graph}

In this section, we \b{relax $3$-connectivity and} study cube-ideal set-systems with connectivity \rt{(}at least\rt{)} $2$\rb{, and}{. More specifically, we study the $2$-cover graph and core of such set-systems, and then} prove \Cref{cube-ideal-size} on the size of their core. \rt{We will need to study the $2$-cover graph and core of such set-systems.}

Take a \emph{preorder} $([n],\succeq)$, i.e., a binary relation $\succeq$ on $[n]$ that satisfies \emph{reflexivity} ($i\succeq i,\,\forall i\in [n]$) and \emph{transitivity} ($i\succeq j, j\succeq k \Rightarrow i\succeq k$, $\forall i,j,k\in [n]$). The \emph{comparability graph} of the preorder is the undirected graph on vertex set $[n]$ with an edge between every pair of distinct vertices that are comparable in the preorder. 

Let $S\subseteq \{0,1\}^n$ be a set-system with connectivity at least $2$. Recall that $G(S)$ denotes the $2$-cover graph of~$S$. Note that $G(S)$ is invariant under twisting $S$. 

\begin{LE}\label{comparability-graph-LE}
Let $S\subseteq \{0,1\}^n$ be a set-system where $\frac{1}{2}\1\in \conv(S)$. Then, after possibly twisting $S$, every $2$-GSC inequality valid for $\conv(S)$ is of the form $x_i\geq x_j$ for distinct $i,j\in [n]$. Furthermore, $G(S)$ is a comparability graph. 
\end{LE}
\begin{proof}
Let us write $\frac{1}{2}\1$ as a convex combination of the points in $S$ with coefficients $\mu_p,p\in S$ where $\mu\geq {\bf 0}$ and $\sum_{p\in S}\mu_p=1$. After twisting the coordinates of $S$, if necessary, we may assume that $\mu_{\bf 0}>0$. Every $2$-GSC inequality valid for $\conv(S)$ is satisfied at equality at $\frac12 \1$, and therefore at ${\bf 0}$, so it must be of the form $x_i\geq x_j$ for some distinct indices $i,j\in [n]$. Consider now the binary relation $\succeq$ on $[n]$ where $i\succeq j$ if $x_i\geq x_j$ is a valid inequality for $\conv(S)$, for $i,j\in [n]$. It can be readily checked that $([n],\succeq)$ is reflexive and transitive, so it is a preorder, and clearly $G(S)$ is its comparability graph, as required.
\end{proof}

We obtain the following fact about the $2$-cover graph of a cube-ideal set-system.

\begin{theorem}\label{comparability-graph}
Let $S\subseteq \{0,1\}^n$ be a cube-ideal set-system with connectivity at least $2$. Then $G(S)$ is a comparability graph.
\end{theorem}
\begin{proof}
As $S$ is cube-ideal with connectivity at least $2$, it follows that $\frac{1}{2}\1\in \conv(S)$. Subsequently, $G(S)$ is a comparability graph by \Cref{comparability-graph-LE}.
\end{proof}

Recall that $\core(S)$ is the set of points in $S$ that satisfy every $2$-GSC inequality valid for $\conv(S)$ at equality. We have the following phenomenon.

\begin{theorem}\label{cube-ideal-core}
Let $S\subseteq \{0,1\}^n$ be a cube-ideal set-system with connectivity at least $2$. Then $\core(S)$ is also cube-ideal.
\end{theorem}
\begin{proof}
Observe that an inequality description of $\conv(\core(S))$ can be obtained from that of $\conv(S)$ by adding, for each $2$-GSC inequality valid for $\conv(S)$, the inequality obtained by flipping the direction of the inequality. For each $2$-GSC inequality, however, the reverse inequality is also GSC: for example, $x_i+x_j\geq 1$ reversed is simply $(1-x_i)+(1-x_j)\geq 1$. Subsequently, $\core(S)$ remains a cube-ideal set-system.
\end{proof}

We now give a description for the convex hull of the core of a cube-ideal set-system $S$. Let $\sum_{i\in I}x_i + \sum_{j\in J}(1-x_j)\geq 1$ be a GSC inequality valid for $\conv(S)$. Recall that the inequality is rainbow if $I\cup J$ intersects every connected component of $G(S)$ at most once. \b{As $G(S)$ is twisting invariant, a GSC inequality being rainbow is also twisting invariant.} As we see below, capacity and rainbow inequalities are sufficient to describe $\conv(\core(S))$.

\begin{LE}\label{rainbow-LE}
Let $S\subseteq \{0,1\}^n$ be a cube-ideal set-system with connectivity at least $2$, and fix an inequality description of $\conv(S)$ comprised of capacity and GSC inequalities. Then the following statements hold: \begin{enumerate}
\item after possibly twisting $S$, we have $x_i=x_j$ for all $x\in \conv(\core(S))$ and indices $i,j$ in the same connected component of $G(S)$,
\item every facet-defining inequality for $\conv(\core(S))$ is equivalent to either a capacity or rainbow inequality in the fixed description of $\conv(S)$.
\end{enumerate}
\end{LE}
\begin{proof}
{\bf (1)} 
By \Cref{comparability-graph-LE}, we may assume after possibly twisting $S$ that every $2$-GSC inequality valid for $\conv(S)$ is of the form $x_i\geq x_j$ for distinct $i,j\in [n]$. In particular, $x_i=x_j$ for all $x\in \conv(\core(S))$ and any pair $i,j$ of indices in the same connected component of $G(S)$. 

{\bf (2)} Take a non-capacity facet-defining inequality $a^\top x\geq \beta$ for $\conv(\core(S))$. A description of the polytope $\conv(\core(S))$ is obtained from the fixed inequality description of $\conv(S)$ after setting all the $2$-GSC inequalities to equality. Subsequently, $a^\top x\geq \beta$ is equivalent to an implicit inequality in the fixed description for $\conv(S)$, say of the form \begin{equation}\label{eq:rainbow-FDI-1}
\sum_{i\in I}x_i + \sum_{j\in J}(1-x_j)\geq 1
\end{equation}
for disjoint subsets $I,J\subseteq [n]$. Denote by $\mathcal{K}$ the set of connected components of $G(S)$, and for each $K\in \mathcal{K}$, take a representative index $i_K\in K$. By part (1), the inequality $\sum_{i\in I}x_i + \sum_{j\in J}(1-x_j)\geq 1$ can be written equivalently as the following valid inequality for {$\conv(\core(S))$}: \begin{equation}\label{eq:rainbow-FDI-2}
\sum_{K\in \mathcal{K}} |I\cap K| x_{i_K} + \sum_{K\in \mathcal{K}} |J\cap K| (1-x_{i_K}) \geq 1.
\end{equation} 
Given that $\1\geq x\geq {\bf 0}$, and the right-hand side value in \eqref{eq:rainbow-FDI-2} is $1$, the nonzero coefficients $|I\cap K|,|J\cap K|$ on the left-hand side can be truncated to $1$ all the while keeping the inequality valid. Since \eqref{eq:rainbow-FDI-2} is facet-defining for $\conv(\core(S))$, and therefore not strictly dominated by another valid inequality, we must have that $|I\cap K|,|J\cap K|\leq 1$ for all $K\in \mathcal{K}$. Furthermore, we cannot have $|I\cap K|=|J\cap K|= 1$ for some $K\in \mathcal{K}$, since otherwise \eqref{eq:rainbow-FDI-2} will be dominated by the equality $x_{i_K} + (1-x_{i_K}) = 1$. Subsequently, $|(I\cup J)\cap K|\leq 1$ for all $K\in \mathcal{K}$, implying in turn that \eqref{eq:rainbow-FDI-1} is a rainbow inequality valid for $\conv(S)$.
\end{proof}

We are now ready to prove \Cref{cube-ideal-size}.

\begin{proof}[Proof of \Cref{cube-ideal-size}]
By \Cref{rainbow-LE} part (1), after possibly twisting $S$, we may assume that $x_i=x_j$ for all $x\in \core(S)$ and indices $i,j$ in the same connected component of $G(S)$. Let $\setcore{S}\in \{0,1\}^d$ be the set-system obtained from $S$ after keeping only one representative from each of the $d$ connected components of $G(S)$, and dropping the remaining indices. Clearly, each rainbow inequality in the description for $\conv(S)$ yields a valid GSC inequality for $\conv(\setcore{S})$ with the same number of variables, and by \Cref{rainbow-LE} part (2), these inequalities together with the capacity inequalities are sufficient to describe $\conv(\setcore{S})$. 
In particular, $S'$ is a cube-ideal set-system with connectivity $\kappa\geq 3$. Thus, $|\core(S)|=|\setcore{S}|\geq 2^{(1-H(1/\kappa))d}$ by \Cref{cube-ideal-size-c3}.
\end{proof}

\section{Ideal clutters}\label{sec:ideal}

In this section, we introduce ideal clutters, connect them to cube-ideal set-systems in two ways, and deduce what our results mean for them; these will be useful in the next section. 

Let $V$ be a finite set of \emph{elements}, and let $\mathcal{C}$ be a family of subsets of $V$ called \emph{members} or \emph{sets}. $\mathcal{C}$ is a \emph{clutter} over \emph{ground set} $V$ if no member contains another~\cite{Edmonds70}. $\mathcal{C}$ is \emph{ideal} if the associated set covering polyhedron, namely $Q(\mathcal{C}):=\left\{x\in \mathbb{R}^V_+:\right.$ $\left.\sum_{v\in C}x_v\geq 1,\,\forall C\in \mathcal{C}\right\}$, has only integral vertices~\cite{Cornuejols94}. 

In integer and linear programming, ideal clutters correspond to set covering linear programs that have an integral optimal solution for any objective vector for which there is a finite optimum. This can be guaranteed when, for instance, the coefficient matrix of the linear system defining $Q(\mathcal{C})$ is \emph{totally unimodular}, i.e., every nonzero sub-determinant is $\pm 1$. In graph and matroid theory, they correspond to multi-commodity flow problems where the \emph{cut condition} necessary for the existence of a flow is also sufficient~\cite{Seymour81,Guenin01,Guenin16}.

To every set-system $S\subseteq \{0,1\}^n$ we can associate a clutter, and through this correspondence we can see a connection between cube-idealness and idealness. 
The \emph{cuboid} of $S$, denoted $\cuboid(S)$, is the clutter over ground set \rb{$[2n]$}{$\{i,\bar{i}:i\in [n]\}$} whose members are \rb{$\{2i-1:p_i=1\}\cup \{2j:p_j=0\}$}{$\{i:p_i=1\}\cup \{\bar{j}:p_j=0\}$}, $\forall p\in S$. It is known $S$ is cube-ideal if, and only if, $\cuboid(S)$ is an ideal clutter~[\cite{Abdi-cuboids}, \b{Theorem 1.6}]. We can use this to prove the following statement claimed in the introduction.

\begin{EG}\label{eg:connectivity-2-linear}
Let $S:=\{{\bf 0},e_1,e_1+e_2,\ldots,e_1+e_2+\cdots+e_n\}\subseteq \{0,1\}^n$. Then $S$ is a cube-ideal set-system with connectivity $2$ such that $\conv(S)$ is full-dimensional and $|S|=n+1$.
\end{EG}
\begin{proof}
It can be readily checked that $\conv(S)$ is full-dimensional, $|S|=n+1$, and that $\frac12 \1$ is the midpoint of the edge of $\conv(S)$ connecting ${\bf 0},\1\in S$. In particular, $S$ must have connectivity $2$. It remains to show that $S$ is cube-ideal. Observe that the incidence matrix of $\cuboid(S)$ has the \emph{consecutive $1$s property}, that is, its columns can be permuted so that the $1$s in each row appear consecutively \b{\cite{Fulkerson64}}. This implies that the coefficient matrix of $Q(\cuboid(S))$ is totally unimodular~\b{[\cite{Schrijver98}, \S19.3, Example 7]}, implying in turn that $\cuboid(S)$ is ideal, thus $S$ is a cube-ideal set-system.
\end{proof}

Let $\mathcal{C}$ be a clutter over ground set $V$. A \emph{cover} is a subset $B\subseteq V$ such that $B\cap C\neq \emptyset$ for all $C\in \mathcal{C}$. The \emph{covering number} of $\mathcal{C}$, denoted $\tau(\mathcal{C})$, is the minimum \rb{cardinality}{size} of a cover. A cover is \emph{minimal} if it does not contain another cover. The \emph{blocker} of $\mathcal{C}$, denoted $b(\mathcal{C})$, is the clutter over ground set $V$ whose members are the minimal covers of $\mathcal{C}$~\cite{Edmonds70}. It is well-known that $b(b(\mathcal{C}))=\mathcal{C}$~\cite{Isbell58,Edmonds70}. Observe that if $\mathcal{C}$ is ideal, then the vertices of $Q(\mathcal{C})$ are precisely the indicator vectors of the minimal covers of $\mathcal{C}$. A fascinating feature of idealness is that it is closed under taking the blocker~\cite{Fulkerson71,Lehman79}. In fact, $\mathcal{C}$ is ideal if, and only if, the \emph{width-length} inequality holds for all $w,\ell\in \mathbb{R}^V_{\geq 0}$, that is, $\min\{w(C):C\in \mathcal{C}\}\cdot \min\{\ell(B):B\in b(\mathcal{C})\}\leq w^\top \ell$~\cite{Lehman79}.

To every clutter we can associate a set-system, and through this correspondence we get yet another connection between idealness and cube-idealness. Let $S(\mathcal{C}):=\{\1_{C}:C\subseteq V, C \text{ contains a set in $\mathcal{C}$}\}\subseteq \{0,1\}^V$. The set-system $S(\mathcal{C})$ is \emph{up-monotone}, that is, if $p\geq q$ for some $p,q\in \{0,1\}^V$ where $q\in S(\mathcal{C})$, then $p\in S(\mathcal{C})$. Note further that the points in $S$ of minimal support are precisely the indicator vectors of the sets in $\mathcal{C}$. It is known that $\mathcal{C}$ is ideal if, and only if, $S(\mathcal{C})$ is a cube-ideal set-system~[\cite{Abdi-cuboids}, \b{Theorem 4.3}]. We have the following.

\begin{theorem}
Let $S\subseteq \{0,1\}^n$ be a cube-ideal set-system with connectivity $\lambda\geq 1$ that is up-monotone. Then $S$ has VC dimension at least $(1-1/\lambda)n$.
\end{theorem}
\begin{proof}
Observe that $S=S(\mathcal{C})$ for a clutter $\mathcal{C}$ over ground set $[n]$. It can be readily checked that $S$ has VC dimension exactly $n-\min\{|C|:C\in \mathcal{C}\}$ and connectivity exactly $\tau(\mathcal{C})=:\lambda$. As $S$ is cube-ideal, it follows that $\mathcal{C}$ is ideal. Thus, by the width-length inequality, $\min\{|C|:C\in \mathcal{C}\} \leq \frac{n}{\min\{\b{|B|}:B\in b(\mathcal{C})\}} = \frac{n}{\lambda}$. Subsequently, $S$ has VC dimension at least $(1-1/\lambda)n$.
\end{proof}

Given an integer $\tau\geq 1$, $\mathcal{C}$ is \emph{$\tau$-cover-minimal} if it has covering number $\tau$, and every element appears in a minimum cover. A clutter is \emph{cover-minimal} if it is $\tau$-cover-minimal for some integer $\tau\geq 1$. The \emph{core} of a cover-minimal clutter $\mathcal{C}$, denoted $\core(\mathcal{C})$, is the clutter of all sets of $\mathcal{C}$ that intersect every minimum cover exactly once. 

\begin{LE}\label{uphull}
Let $\mathcal{C}$ be an ideal $\tau$-cover-minimal clutter \b{over ground set $V$}, for some integer $\tau\geq 1$. Then $S(\mathcal{C})$ is a cube-ideal set-system with connectivity $\tau$\rt{, where every variable appears in a valid $\tau$-SC inequality}. Furthermore, for $F$ the minimal face of $\conv(S(\mathcal{C}))$ containing $\frac{1}{\tau} \1$, the dimension of $F$ is $|V|-r$ where $r$ is the rank of $\{\1_B:B\in b(\mathcal{C}), |B|=\tau\}$, and $S(\mathcal{C})\cap F = \{\1_C:C\in \core(\mathcal{C})\}$.
\end{LE}
\begin{proof}
Let $S:=S(\mathcal{C})$. As $\mathcal{C}$ is ideal, $S$ is a cube-ideal set-system whose convex hull is described by $\0\leq x\leq \1, x(B)\geq 1,\,\forall B\in b(\mathcal{C})$. In particular, $S$ has connectivity $\tau(\mathcal{C})=\tau$\rt{, and as $\mathcal{C}$ is $\tau$-cover-minimal, every variable appears in a $\tau$-SC inequality valid for $S$}. Furthermore, the minimal face of $\conv(S)$ containing $\frac{1}{\tau} \1$ is \begin{align*}
F
&=
\b{[0,1]^V \cap 
\left\{x\in \mathbb{R}^V:x(B)\geq 1,\forall B\in b(\mathcal{C}) \text{ s.t. } |B|>\tau\right\}}\cap 
\left\{x\b{\in \mathbb{R}^V}:x(B)=1,\forall B\in b(\mathcal{C}) \text{ s.t. } |B|=\tau\right\}
\end{align*} where {the third set contains all the implicit equalities}. 
Subsequently, the dimension of $F$ is $|V|-r$. It remains to prove $S\cap F = \{\1_C:C\in \core(\mathcal{C})\}$. The inclusion $\supseteq$ is clear. For the reverse inclusion, suppose $\widehat{C}\subseteq V$ contains a set $C\in \mathcal{C}$, and $\1_{\widehat{C}}\in S\cap F$. We claim that $\widehat{C}=C$. For if not, pick $v\in \widehat{C}\setminus C$, and let $B$ be a minimum cover of $\mathcal{C}$ that contains $v$, which exists as $\mathcal{C}$ is cover-minimal. As $B\cap C\neq \emptyset$ and $v\in B$, it follows that $|B\cap \widehat{C}|\geq 2$, a contradiction as $|B\cap \widehat{C}|=1$.
\end{proof}

We may therefore apply \Cref{cube-ideal-faces-lower-bound} to obtain the following inequality. 

\begin{theorem}\label{ideal-core-exponential}
For every integer $\tau\geq 3$ and $\beta\in (0,1)$, there is a constant $\theta\b{:=\theta_{\ref{cube-ideal-faces-lower-bound}}(\tau,\beta)}>0$ such that the following statement holds: 

Let $\mathcal{C}$ be an ideal $\tau$-cover-minimal clutter over ground set $V$. 
Suppose the rank of $\{\1_B:B\in b(\mathcal{C}), |B|=\tau\}$ is at most $(1-\beta)|V|$. Then $|\core(\mathcal{C})|\geq e^{\theta |V|}$.
\end{theorem}
\begin{proof}
Let $S:=S(\mathcal{C})$. By \Cref{uphull}, 
$S$ is a cube-ideal set-system with connectivity $\lambda:=\tau$; furthermore, for $F$ the minimal face of $\conv(S)$ containing $\frac{1}{\lambda} \1$, the dimension of $F$ is at least $|V| - (1-\beta)|V| = \beta|V|$, and $|\core(\mathcal{C})| = |S\cap F|$. Thus, by \Cref{cube-ideal-faces-lower-bound}, $|\core(\mathcal{C})| \geq e^{\theta |V|}$ for $\theta = \b{\theta_{\ref{cube-ideal-faces-lower-bound}}}(\lambda,\beta)$.
\end{proof}

Let $\mathcal{C}$ be a $2$-cover-minimal clutter over ground set $V$, these are also known as \emph{tangled} clutters~\cite{Abdi-kwise}. 
\b{As a consequence of the tools developed so far, we prove that $\core(\mathcal{C})$ is ideal and give an exponential lower bound on $|\core(\mathcal{C})|$; let us elaborate.} \r{[line break]}

Denote by $G(\mathcal{C})$ the graph over vertex set $V$ whose edges correspond to the minimum covers of $\mathcal{C}$. A cover of $\mathcal{C}$ is \emph{rainbow} if it intersects every connected component of $G(\mathcal{C})$ at most once. 
The \emph{rainbow covering number of $\mathcal{C}$}, denoted $\mu({\mathcal{C}})$, is the minimum size of a rainbow cover of $\mathcal{C}$; if there is no rainbow cover, then $\mu({\mathcal{C}}):=+\infty$. Observe that $\mu({\mathcal{C}})\geq 3$. \rt{We prove the following exponential lower bound on the size of the core of such clutters.}

\b{Given disjoint subsets $I,J\subseteq V$, the \emph{minor obtained after deleting $I$ and contracting $J$}, denoted by $\mathcal{C}\setminus I/J$, is the clutter over ground set $V\setminus (I\cup J)$ whose sets are the inclusion-wise minimal sets in $\{C\setminus J:C\in \mathcal{C},C\cap I=\emptyset\}$. If $\mathcal{C}$ is an ideal clutter, then so are its minors~(see \cite{Cornuejols01}, \S1.4).}

\begin{theorem}\label{main}
Let $\mathcal{C}$ be an ideal $2$-cover-minimal clutter \b{over ground set $V$}, let $d$ be the number of connected components of $G(\mathcal{C})$, and let $\mu$ be its rainbow covering number. 
Then \b{$\core(\mathcal{C})$ is an ideal clutter such that} $|\core(\mathcal{C})|\geq 2^{(1-H(1/\mu))d}$.
\end{theorem}
\begin{proof}
Let $S:=S(\mathcal{C})$. We know from \Cref{uphull} that $S$ is a cube-ideal set-system with connectivity $2$\rb{. F}{; f}urthermore, every valid $2$-GSC inequality \rb{is in fact a $2$-SC inequality}{for $S$ is of the form $x(B)\geq 1$ for a minimum cover $B$ of $\mathcal{C}$}, so $G(S)=G(\mathcal{C})$ and $\core(S) = \{\1_C:C\in \core(\mathcal{C})\}$. The claimed inequality \rt{now} follows from \Cref{cube-ideal-size}. 
\b{For the first part of the claim, note that $\core(S)$ is a cube-ideal set-system by \Cref{cube-ideal-core}, so $\cuboid(\core(S))$ is an ideal clutter~[\cite{Abdi-cuboids}, Theorem 1.6]. As $\core(\mathcal{C})$ is a minor of $\cuboid(\core(S))$ obtained after contracting $\{\bar{v}:v\in V\}$, it is also an ideal clutter.} 
\end{proof}

In a similar fashion, \Cref{cube-ideal-face-CN} implies the following\rt{(the two conjectures are in fact equivalent)}.

\begin{CN}\label{tau-cover-min-CN}
For every integer $\tau\geq 3$, there is a constant $\theta\b{:=\theta_{\ref{cube-ideal-face-CN}}(\tau)}>0$ such that the following statement holds: 

Let $\mathcal{C}$ be an ideal $\tau$-cover-minimal clutter over ground set $V$.  
Let $d$ be the dimension of the minimal face of $Q(b(\mathcal{C}))$ containing $\frac{1}{\tau} \1$. 
Then $|\core(\mathcal{C})|\geq e^{\theta d}$.
\end{CN}

\section{Applications to combinatorial optimization}\label{sec:apps}

In this section, we present three further applications of our results to combinatorial optimization.

\subsection{Strong orientations of mixed graphs}\label{subsec:mixed-graph}

Here we prove \Cref{strong-orientations-exp}. We first show that strong orientations of a mixed graph correspond to a cube-ideal set-system. Let us elaborate.

Let $\mathcal{F}$ be a family of subsets of a finite set $V$. $\mathcal{F}$ is a \emph{crossing family} if $U\cap W,U\cup W\in \mathcal{F}$ for any two sets $U,W\in \mathcal{F}$ such that $U\cap W\neq \emptyset$ and $U\cup W\neq V$. For a crossing family $\mathcal{F}$, a function $f:\mathcal{F}\to \mathbb{Z}$ is \emph{crossing supermodular} if $f(U\cap W)+f(U\cup W)\geq f(U)+f(W)$ for any two sets $U,W\in \mathcal{F}$ such that $U\cap W\neq \emptyset$ and $U\cup W\neq V$.

Let $A\in\mathbb{Q}^{m\times n}$ and $b\in \mathbb{Q}^m$. The system $Ax\leq b$ is \emph{totally dual integral (TDI)} if for every cost vector $c\in \mathbb{Z}^n$ for which $\max\{c^\top x: Ax\leq b\}$ has an optimum solution, its dual $\min\{ b^\top y: A^\top y = c, y\geq 0\}$ has an integral optimum solution. If the right-hand side $b$ is integral and the system $Ax\leq b$ is TDI, then the polyhedron $P:= \{x: Ax\leq b\}$ is integral~\cite{Edmonds77}. A system $Ax\leq b$ is called \emph{box-TDI} if for every pair of \rb{integral vectors}{vectors $\ell,u\in \mathbb{Z}^n\cup \{\pm \infty\}^n$ such that} $\ell\leq u$, the system $Ax\leq b, \ell\leq x\leq u$ is TDI. In particular, every box-TDI system is also TDI. 

Let $G=(V,A,E)$ be a mixed graph where every pseudo dicut contains at least two edges from~$E$. Let $\overrightarrow{E}$ be an arbitrary reference orientation of $E$. \rb{Denote by}{Recall that} $\scr(G;\overrightarrow{E})$ \b{is} the set of vectors $x\in \{0,1\}^{\overrightarrow{E}}$ such that the re-orientation of $\overrightarrow{E}$ obtained after flipping the arcs in $\{a\in \overrightarrow{E}:x_a=1\}$ is $G$-strong, i.e., turns $G$ into a strongly connected digraph.

\begin{theorem}\label{SCR-cube-ideal}
For a mixed graph $G=(V,A,E)$, and a reference orientation $\overrightarrow{E}$, the convex hull of $\scr(G;\overrightarrow{E})$ is described by 
\begin{align}
\sum_{a\in \delta_{\overrightarrow{E}}^+(U)} x_a + \sum_{b\in \delta_{\overrightarrow{E}}^-(U)} (1-x_b)&\geq 1 \quad \forall U\subset V, U\neq \emptyset, \delta_A^-(U)=\emptyset,\label{eq:supflow-1}\\
\1\geq x&\geq {\bf 0}.\label{eq:supflow-2}
\end{align} In particular, the set-system $\scr(G;\overrightarrow{E})$ is cube-ideal. 
\end{theorem}
\begin{proof}
\rb{It can be readily checked that t}{T}he integral solutions to \eqref{eq:supflow-1}-\eqref{eq:supflow-2} are precisely the points in $\scr(G;\overrightarrow{E})$, \b{since $x\in\{0,1\}^{\overrightarrow{E}}$ satisfies \eqref{eq:supflow-1} iff the corresponding re-orientation is $G$-strong}. Thus, it remains to show that the system is integral. To this end, let $\mathcal{F}:=\{U\subset V: U\neq \emptyset, \delta_A^-(U)=\emptyset\}$, which is a crossing family. Note that \eqref{eq:supflow-1} can be rewritten as \begin{equation}\label{eq:supflow-3}
x(\delta_{\overrightarrow{E}}^+(U))-x(\delta_{\overrightarrow{E}}^-(U))\geq 1-|\delta_{\overrightarrow{E}}^-(U)| \quad\forall U\in \mathcal{F}.
\end{equation} Given that $U\in \mathcal{F}\mapsto 1-|\delta_{\overrightarrow{E}}^-(U)|\in \mathbb{Z}$ is a crossing supermodular function, we conclude from a well-known result~\cite{Edmonds77} that \eqref{eq:supflow-3} is box-TDI. In particular, \eqref{eq:supflow-1}-\eqref{eq:supflow-2} is TDI, and is therefore integral, as required.
\end{proof}

Suppose every pseudo dicut of $G$ contains at least two edges from $E$. The \emph{$2$-pseudo-dicut graph of $G$} is the graph~$H$ on vertex set $E$ with an edge between distinct $e,f$ for every pseudo dicut of $G$ that contains only $e$ \b{and} $f$ from~$E$. A pseudo dicut of $G$ is \emph{rainbow} if it contains at most one element of $E$ from every connected component of the $2$-pseudo-dicut graph. In particular, every rainbow pseudo dicut contains at least $3$ elements from $E$. We are now ready to prove the following theorem. Note that $\kappa:=+\infty$ if there is no rainbow pseudo dicut.

\begin{theorem}\label{strong-orientations}
Let $G=(V,A,E)$ be a mixed graph where every pseudo dicut contains at least $2$ edges from~$E$. Let $d$ be the number of connected components of the $2$-pseudo-dicut graph of $G$, and let $\kappa$ be the minimum number of elements from $E$ in a rainbow pseudo dicut. Then the number of $G$-strong orientations is at least $2^{(1-H(1/\kappa))d}$.
\end{theorem}
\begin{proof}
Let $\overrightarrow{E}$ be an arbitrary reference orientation of $E$, and let $S:=\scr(G;\overrightarrow{E})$, which is cube-ideal by \Cref{SCR-cube-ideal}. Fix the inequality description \eqref{eq:supflow-1}-\eqref{eq:supflow-2} for $\conv(S)$. It therefore follows from \Cref{cube-ideal-size} that $|\core(S)|\geq 2^{(1-H(1/\kappa'))d'}$, where $d'$ is the number of connected components of $G(S)$, and $\kappa'$ is the minimum number of variables used in a rainbow inequality in \eqref{eq:supflow-1}. Observe that $G(S)$ is precisely the $2$-pseudo-dicut graph of $G$. Furthermore, the rainbow inequalities in \eqref{eq:supflow-1} correspond to the rainbow pseudo dicuts of $G$, and the number of variables in the inequality is equal to the number of edges from $E$ in the pseudo dicut. Subsequently, $d'=d$ and $\kappa'=\kappa$, so $|S|\geq |\core(S)|\geq 2^{(1-H(1/\kappa))d}$, as required.
\end{proof}

We obtain \Cref{strong-orientations-exp} as a consequence of \rb{this theorem}{\Cref{strong-orientations}}.

\begin{proof}[Proof of \Cref{strong-orientations-exp}]
Let $G=(V,A,E)$ be a mixed graph where every pseudo dicut contains at least $\lambda$ edges from~$E$, where $\lambda\geq 3$. Then the $2$-pseudo-dicut graph of $G$ has no edge, so it has exactly $|E|$ connected components, and every pseudo dicut of $G$ is rainbow. Subsequently, the number of $G$-strong orientations is at least $2^{(1-H(1/\lambda))|E|}$ \b{by \Cref{strong-orientations}}.
\end{proof}

\subsection{Dijoins in $\tau$-sink-regular bipartite digraphs}\label{subsec:dijoin}

\begin{proof}[Proof of \Cref{tight-dijoins-exponential}]
Let $\mathcal{C}$ be the clutter of minimal dijoins of $D$. It is known that $\mathcal{C}$ is an ideal clutter, and that $b(\mathcal{C})$ consists of the minimal dicuts of $D$~[\cite{Lucchesi78}, \b{see also \cite{Cornuejols01}, \S1.3.4}]. Thus, $\mathcal{C}$ is $\tau$-cover-minimal \b{as $D$ is a bipartite digraph where every sink has degree $\tau$ and every dicut has size at least $\tau$}, and $\core(\mathcal{C})$ is precisely the clutter of minimal dijoins of $D$ intersecting every minimum dicut exactly once. We shall therefore apply \Cref{ideal-core-exponential} to argue that $|\core(\mathcal{C})|$ is exponentially large. To this end, let $\mathcal{F}:=\{U\subset V:U\neq \emptyset, |\delta^+(U)|=\tau, \delta^-(U)=\emptyset\}$, and 
let $M$ be the matrix whose rows correspond to $\{\1_{\delta^+(U)}:U\in \mathcal{F}\}$. We claim that $M$ has row rank at most $(1 - \beta) |A|$ for $\beta=\frac13$, thus allowing us to apply \Cref{ideal-core-exponential} \b{for $\theta:=\theta_{\ref{cube-ideal-faces-lower-bound}}(\tau,\frac13)$} to finish the proof. Let $S,T$ be the sets of sources and sinks of $D$, respectively. Then we have $$
\1_{\delta^+(U\cap S)} - \1_{\delta^-(U\cap T)} = \1_{\delta^+(U)} - \1_{\delta^-(U)} = \1_{\delta^+(U)} \qquad \forall U\in \mathcal{F}.
$$ Note that $\1_{\delta^-(v)},\forall v\in T$ are rows of $M$. Thus, after applying elementary row operations to the rows of $M$ corresponding to $\1_{\delta^+(U)},U\in \mathcal{F}, U\not\supseteq S$, we obtain a matrix $N$ whose rows are 
$\1_{\delta^-(v)},\,\forall v\in T;\,\1_{\delta^+(U\cap S)},\,\forall U\in \mathcal{F}, U\not\supseteq S.$ The first set of rows of $N$ clearly has rank at most $|T|${. The second set of rows is generated by $\1_{\delta^+(u)},u\in S$, and therefore has rank at most $|S|$.}
\rb{implying in turn that}{Consequently,} $N$ and therefore $M$ has rank at most $|S|+|T|=|V|$. Given that every source has degree at least $\tau$ (as the arcs incident with it form a dicut), it follows that $|S|\leq |T|$, so $|V|\leq 2|T| = \frac{2|A|}{\tau}\leq \frac23 |A|$, thus proving the claim. 
\end{proof}

\subsection{Perfect matchings in $r$-graphs \b{and the Lov\'{a}sz--Plummer conjecture}}\label{subsec:perfect-matching}

\Cref{ideal-core-exponential} has a consequence for the number of perfect matchings in certain graphs. To elaborate, let $r\geq 3$ be an integer. An \emph{$r$-graph} is an $r$-regular graph $G=(V,E)$ where $|V|$ is even, and every \emph{odd cut} has size at least~$r$, that is, $|\delta(U)|\geq r$ for all $U\subset V$ where $|U|$ is odd. It is known that every $r$-graph is \emph{matching-covered}, that is, every edge belongs to a perfect matching~\cite{Seymour79}. \b{A well-known conjecture by Lov\'{a}sz and Plummer states that for every $r\geq 3$, there exist constants $c_1(r)>1,c_2(r)>0$ such that every $r$-regular matching-covered graph $G=(V,E)$ contains at least $c_2(r) \cdot c_1(r)^{|V|}$ perfect matchings; furthermore, $c_1(r)\to \infty$ as $r\to \infty$~[\cite{Lovasz86}, Conjecture 8.1.8].
The first part of this conjecture was proved for $3$-graphs~\cite{Esperet11}. In this subsection, we prove (a refinement of) this part for all $r$-graphs for $r\geq 4$.}

Perfect matchings intersecting every minimum odd cut exactly once correspond to the core of a certain ideal $r$-cover-minimal clutter. Through this connection, we shall prove the following as a consequence of \Cref{ideal-core-exponential}.

\begin{theorem}\label{perfect-matching-exp}
For every integer $r\geq 3$ and $\beta\in (0,1)$, there is a constant $\theta\b{:=\theta_{\ref{cube-ideal-faces-lower-bound}}(r,\beta)}>0$ such that the following statement holds: 

Let $G=(V,E)$ be an $r$-graph where the rank of $\{\1_{\delta(U)}:|\delta(U)|=r,\, |U| \text{ is odd}\}$ is at most $(1 - \beta) |E|$. Then the number of perfect matchings of $G$ that intersect every minimum odd cut exactly once is at least $e^{\theta |E|}$.
\end{theorem}
\begin{proof}
A \emph{postman set} is a subset $J\subseteq E$ where every vertex of $G$ is incident to an odd number of edges from $J$. Let $\mathcal{C}$ be the clutter over ground set $E$ of the minimal postman sets of $G$. It is known that $\mathcal{C}$ is an ideal clutter, and that $b(\mathcal{C})$ consists of the minimal odd cuts of $G$~\cite{Edmonds73}. In particular, as $G$ is an $r$-graph, $\mathcal{C}$ is an $r$-cover-minimal clutter, and $\core(\mathcal{C})$ is precisely the clutter of perfect matchings of $G$ that intersect every minimum odd cut exactly once. The claim now follows from \Cref{ideal-core-exponential} \b{for $\theta:=\theta_{\ref{cube-ideal-faces-lower-bound}}(r,\beta)$}.
\end{proof}

\b{By applying \Cref{barvinok} to a different polytope,} Barvinok showed that every $3$-graph $G=(V,E)$ that is \emph{essentially $4$-edge-connected}, i.e., $|\delta(U)|\geq 4$ for all $U\subset V$ such that $1<|U|<|V|-1$, has exponentially many perfect matchings~\cite{Barvinok13}. \b{Observe that in this case, the rank of $\{\1_{\delta(U)}:|\delta(U)|=3,\, |U| \text{ is odd}\}$ is at most $|V| = \frac{2}{3}|E|$, so we can also apply \Cref{perfect-matching-exp} for $r=3,\beta=\frac13$ to obtain this result.

Using graph theoretic techniques,} Esperet et al.\ showed that the number of perfect matchings in \emph{any} $3$-graph $G=(V,E)$ is at least $2^{|V|/3656}$. The authors then used this result to show that every $(r-1)$-edge-connected $r$-graph for $r\geq 4$, has at least $2^{f(r)\cdot |V|}$ perfect matchings for $f(r) = \frac{1}{3656}\left(1-\frac{1}{r}\right)\left(1-\frac{2}{r}\right)$~\cite{Esperet11}.

\b{
It turns out that for $r\geq 4$, \emph{all} $r$-graphs satisfy the rank condition in \Cref{perfect-matching-exp}. To see this, we need a few notions and lemmas. A family $\mathcal{L}$ over ground set $V$ is \emph{laminar} if for all $U,W\in \mathcal{L}$ such that $U\cap W\neq \emptyset$, we have $U\subseteq W$ or $W\subseteq U$. A simple inductive argument implies that a laminar family $\mathcal{L}$ has size at most $2|V|-1$. However, this upper bound can be improved in the following relevant case.}

\b{
\begin{LE}\label{laminar-odd}
Let $\mathcal{L}$ be a laminar family of odd subsets of $\{1,\ldots,2n\}$. Then $|\mathcal{L}|\leq 3n-1$.
\end{LE}
\begin{proof}
We prove this by induction on $n$. The base case $n=1$ is clear. For the induction step, assume that $n\geq 2$.

Suppose in the first case that some maximal set $A$ in $\mathcal{L}$ has size $\leq 2n-3$. We may clearly assume that $|A|\geq 3$. Let $B:=\overline{A}$. Then $|B|$ is odd, and $\mathcal{L}\cup \{B\}$ is also laminar, so we may assume that $B\in \mathcal{L}$. 
Let 
$\mathcal{L}_1,\mathcal{L}_2$ be obtained from $\mathcal{L}$ after shrinking $B,A$ to a single element $b,a$, respectively. Observe that $\mathcal{L}_1$ contains both $A$ and $\{b\}$, and $\mathcal{L}_2$ contains both $B$ and $\{a\}$. Given that the ground set of each $\mathcal{L}_i$ has $2n_i<2n$ elements, we can apply the induction hypothesis to conclude that $$
|\mathcal{L}| = |\mathcal{L}_1|+|\mathcal{L}_2| -2 \leq (3n_1-1)+(3n_2-1)-2 = 3(n_1+n_2)-4 = 3(n+1)-4=3n-1,
$$ thus completing the induction step in this case.

Suppose in the remaining case that every maximal set in $\mathcal{L}$ has size $2n-1$. Thus, there is a unique maximal set $A$ in $\mathcal{L}$. Now, let $B$ be a proper maximal subset of $A$ in $\mathcal{L}$, which we may assume exists and satisfies $|B|\geq 3$. By definition, $|B|\leq |A|-2$. Then $\mathcal{L}':=(\mathcal{L}\setminus \{A\}) \cup \{\overline{B}\}$ is also laminar family of odd subsets of $\{1,\ldots,2n\}$. Furthermore, $B$ is a maximal set in $\mathcal{L}'$ of size $\leq 2n-3$. Therefore, by falling back to the previous case, we obtain that $|\mathcal{L}|=|\mathcal{L}'|\leq 3n-1$, thus completing the induction step.
\end{proof}

The \emph{perfect matching polytope} of a graph $G=(V,E)$ on an even number of vertices is the convex hull of the indicator vectors of its perfect matchings. This polytope can be described as the set of all vectors $x\in \mathbb{R}^E_{\geq 0}$ such that $x(\delta(v))=1$ for every vertex $v$, and $x(\delta(U))\geq 1$ for every odd cut $\delta(U)$~\cite{Edmonds65}.

\begin{LE}\label{rankUB}
Let $G=(V,E)$ be a matching-covered graph, and let $x^\star$ be a point in the perfect matching polytope where $x^\star_e>0$ for all $e\in E$. Let $
\mathrm{rk}$ be the rank of $\{\1_{\delta(U)} : x^\star(\delta(U))=1, |U| \text{ is odd}\}$. Then $\mathrm{rk}\leq \frac32 |V| - 1$.
\end{LE}
\begin{proof}
Let $\mathcal{F}:=\{\{U,\overline{U}\} : x^\star(\delta(U))=1, |U| \text{ is odd}\}$ be the family of \emph{$x^\star$-tight cuts}. Then $\mathcal{F}$ is a crossing family, in the sense that for all pairs $\{U,\overline{U}\},\{W,\overline{W}\}\in \mathcal{F}$ that \emph{cross}, i.e.\ $U\cap W,U\cap \overline{W},\overline{U}\cap W,\overline{U}\cap \overline{W}\neq \emptyset$, and $|U\cap W|$ is odd, we have that $\{U\cap W,\overline{U\cap W}\},\{U\cup W,\overline{U\cup W}\}\in \mathcal{F}$. To see this, note that $$
2 = x^\star(\delta(U))+x^\star(\delta(W))\geq x^\star(\delta(U\cap W)) + x^\star(\delta(U\cup W))\geq 2,
$$ where the first (in)equality follows from $\{U,\overline{U}\},\{W,\overline{W}\}\in \mathcal{F}$, the second from submodularity of the cut function, and the third from the fact that $|U\cap W|$, and therefore $|U\cup W|$, is odd, and $x^\star$ is in the perfect matching polytope. Subsequently, equality must hold throughout, so $\{U\cap W,\overline{U\cap W}\},\{U\cup W,\overline{U\cup W}\}\in \mathcal{F}$. Observe further that as $x^\star_e>0$ for all $e\in E$, we must have that $
\1_{\delta(U)} + \1_{\delta(W)} = \1_{\delta(U\cap W)}+\1_{\delta(U\cup W)}$. 

Now let $\mathcal{C}$ be a maximal subset of $\mathcal{F}$ that is \emph{cross-free}, i.e. contains no two pairs that cross. 
Note that $\mathcal{C}$ contains all \emph{trivial cuts}, i.e.\ $\delta(v),v\in V$. 
Let $L$ be the linear hull of $\1_{\delta(U)}, \{U,\overline{U}\}\in \mathcal{F}$. 

\begin{claim*} 
$\mathcal{C}$ generates $L$, that is, $L$ is the linear hull of $\1_{\delta(U)}, \{U,\overline{U}\}\in \mathcal{C}$.
\end{claim*}
\begin{cproof}
We prove this by induction on $|V|$. The base case $|V|\leq 4$ is clear as $\mathcal{C}=\mathcal{F}$. For the induction step, assume that $|V|\geq 6$.

Pick an arbitrary nontrivial cut $\{U,\overline{U}\}$ in $\mathcal{C}$. 
Let $G_1,G_2$ be obtained from $G$ after shrinking $U,\overline{U}$ to a single vertex, respectively. 
Let $x^\star_i$ be the restriction of $x^\star$ to $E(G_i)$, which also belongs to the perfect matching polytope of $G_i$, let $\mathcal{F}_i$ be the family of $x^\star_i$-tight cuts, and let $L_i$ be the corresponding linear hull. Observe that $\mathcal{F}_1$ consists of those cuts $\{R,\overline{R}\}$ in $\mathcal{F}$ such that $U\subseteq R$, and $\mathcal{F}_2$ of those cuts $\{R,\overline{R}\}$ in $\mathcal{F}$ such that $\overline{U}\subseteq R$. Let $\mathcal{C}_1$ be the family of cuts $\{R,\overline{R}\}$ in $\mathcal{C}$ such that $U\subseteq R$, and $\mathcal{C}_2$ of those cuts $\{R,\overline{R}\}$ in $\mathcal{C}$ such that $\overline{U}\subseteq R$. The maximality of $\mathcal{C}$ implies that each $\mathcal{C}_i$ is a maximal cross-free subset of $\mathcal{F}_i$. Thus, by the induction hypothesis, each $\mathcal{C}_i$ generates $L_i$.

Now, to prove that $\mathcal{C}$ generates $L$, pick an arbitrary $\{W,\overline{W}\}$ in $\mathcal{F}$ outside $\mathcal{C}$. By the maximality of $\mathcal{C}$, $\{W,\overline{W}\}$ must cross some nontrivial cut $\{U,\overline{U}\}\in \mathcal{C}$, as above. Thus, $\{U\cap W,\overline{U\cap W}\},\{U\cup W,\overline{U\cup W}\}$ also belong to $\mathcal{F}$, and $$
\1_{\delta(U)} + \1_{\delta(W)} = \1_{\delta(U\cap W)}+\1_{\delta(U\cup W)}.
$$ Note that $\delta(U\cap W)$ is an $x^\star_2$-tight cut from $\mathcal{F}_2$, so $\1_{\delta(U\cap W)}\in L_2$, and similarly $\1_{\delta(U\cup W)}\in L_1$. By the induction hypothesis, $\mathcal{C}_2$ generates $\1_{\delta(U\cap W)}$, and $\mathcal{C}_1$ generates $\1_{\delta(U\cup W)}$. In particular, as $\mathcal{C}_i\subseteq \mathcal{C}$, it follows that $\1_{\delta(U\cap W)},\1_{\delta(U\cup W)}$ are generated by $\mathcal{C}$, so the vector equality above implies that $\1_{\delta(W)}$ is also generated by $\mathcal{C}$, as required.
\end{cproof}

In particular, $\mathrm{rk}\leq |\mathcal{C}|$. By carefully retaining exactly one shore from each pair $\{U,\overline{U}\}$ in the cross-free family $\mathcal{C}$, we obtain a laminar family $\mathcal{L}$ of odd subsets of $V$. Subsequently, by \Cref{laminar-odd}, $\mathrm{rk}\leq |\mathcal{C}|=|\mathcal{L}|\leq \frac{3}{2} |V|-1$, thus finishing the proof.
\end{proof}

We are now ready to provide the following exponential lower bound for perfect matchings of an $r$-graph with $r\geq 4$.

\begin{CO}\label{perfect-matching-exp->=4-graph}
For every integer $r\geq 4$, there is a constant $\theta:=\theta_{\ref{cube-ideal-faces-lower-bound}}(r,1-\frac{3}{r})>0$ such that the following statement holds: 

The number of perfect matchings of an $r$-graph $G=(V,E)$ that intersect every minimum odd cut exactly once is at least $e^{\theta |E|}$.
\end{CO}
\begin{proof}
Let $x^\star=\frac{1}{r}\1$, and let $\mathrm{rk}$ be the rank of the minimum odd cuts of $G$. It follows from \Cref{rankUB} that $\mathrm{rk}\leq \frac32 |V| - 1 = \frac{3}{r}|E|-1$ since $G$ is $r$-regular. The theorem now follows from \Cref{perfect-matching-exp} for $\beta:=1-\frac{3}{r}>0$ and $\theta:=\theta_{\ref{cube-ideal-faces-lower-bound}}(r,\beta)>0$.
\end{proof}
}

\begin{figure}
\centering
\includegraphics[scale=0.3]{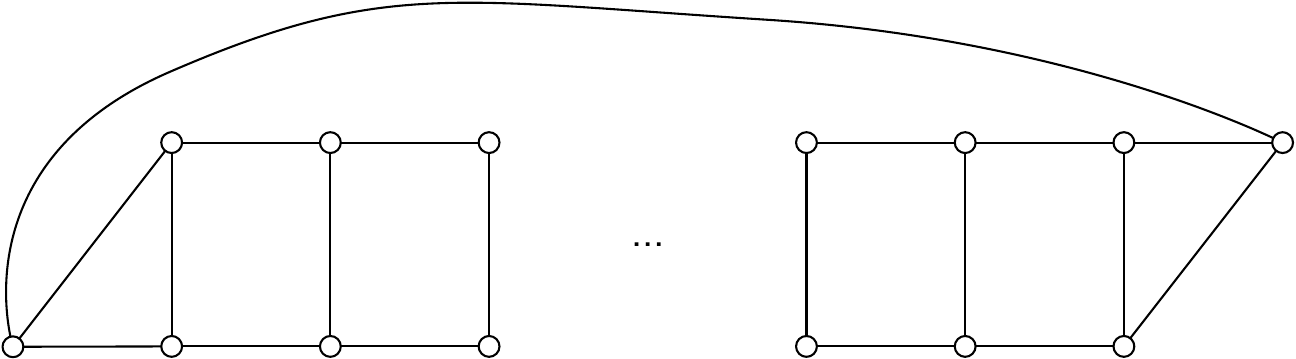}
\caption{A staircase on at least $4$ vertices. This graph has exactly $3$ perfect matchings that intersect every minimum odd cut exactly once.}
\label{fig:staircase}
\end{figure}

\b{\Cref{perfect-matching-exp->=4-graph} does not naturally generalize to $r=3$, as} there exists an unbounded family of $3$-graphs with exactly $3$ perfect matchings that intersect every minimum odd cut exactly once, namely \emph{staircases} as depicted in \Cref{fig:staircase}. This family also shows that the rank condition in \Cref{perfect-matching-exp} cannot be dropped; for such graphs, the rank of $\{\1_{\delta(U)}:|\delta(U)|=3,\, |U| \text{ is odd}\}$ is $|E| - 2$. 

\b{\Cref{tau-cover-min-CN} tells us how we might be able to extend \Cref{perfect-matching-exp->=4-graph} to $r=3$, as follows.

\begin{CN}
The number of perfect matchings of a $3$-graph $G=(V,E)$ that intersect every minimum odd cut exactly once is at least $e^{\theta (|E|-\mathrm{rk})}$, where $\theta=\theta_{\ref{cube-ideal-face-CN}}(3)>0$ and $\mathrm{rk}$ is the rank of $\{\1_{\delta(U)}:|\delta(U)|=3,\, |U| \text{ is odd}\}$.
\end{CN}
}

\section*{Acknowledgements}

Parts of this material are based upon work supported by the National Science Foundation under Grant No. DMS-1929284 while AA, GC and DD were in residence at and/or visiting the Institute for Computational and Experimental Research in Mathematics in Providence, RI, during the Discrete Optimization: Mathematics,  Algorithms, and Computation semester program. This work was also supported by ONR grant 14-22-1-2528 and EPSRC grant EP/X030989/1. We would like to thank Louis Esperet, Stefan Weltge, \b{and Giacomo Zambelli} for fruitful discussions about parts of this work. \b{Finally, we would like to thank anonymous reviewers for their detailed and constructive feedback on an early draft of this manuscript.}

\paragraph{Data Availability Statement.} No data are associated with this article. Data sharing is not applicable to this article.

{\small 
\bibliographystyle{alpha}
\bibliography{cube-ideal-size-references}}

@article{Merino99,
	author = {Merino, C. and Welsh, D. J. A.},
	da = {1999/06/01},
	doi = {10.1007/BF01608795},
	id = {Merino1999},
	isbn = {0219-3094},
	journal = {Annals of Combinatorics},
	number = {2},
	pages = {417--429},
	title = {Forests, colorings and acyclic orientations of the square lattice},
	ty = {JOUR},
	url = {https://doi.org/10.1007/BF01608795},
	volume = {3},
	year = {1999}}

@book{EM22,
	editor = {Ellis-Monaghan, J. A. and Moffatt, I.},
	publisher = {Chapman and Hall/CRC Press},
	title = {Handbook of the Tutte Polynomial and Related Topics},
	year = {2022}}

@article{Jaeger90,
	author = {Jaeger, F. and Vertigan, D. L. and Welsh, D. J. A.},
	booktitle = {Mathematical Proceedings of the Cambridge Philosophical Society},
	db = {Cambridge Core},
	doi = {DOI: 10.1017/S0305004100068936},
	dp = {Cambridge University Press},
	et = {2008/10/24},
	isbn = {0305-0041},
	number = {1},
	pages = {35-53},
	publisher = {Cambridge University Press},
	title = {On the computational complexity of the Jones and Tutte polynomials},
	ty = {JOUR},
	url = {https://www.cambridge.org/core/product/0EA052341269A2816C36B15380B8AA02},
	volume = {108},
	year = {1990}}

@article{Edmonds65,
	author = {Edmonds, J.},
	journal = {Journal of Research of the National Bureau of Standards B},
	pages = {125-130},
	title = {Maximum matching and a polyhedron with 0, 1-vertices},
	volume = {69},
	year = {1965}}

@book{Schrijver98,
	author = {Alexander Schrijver},
	month = {April},
	publisher = {Wiley},
	title = {Theory of Linear and Integer Programming},
	year = {1998}}

@article{Whitney32,
	author = {Whitney, H.},
	journal = {American Journal of Mathematics},
	pages = {150--168},
	title = {Non-separable and planar graphs},
	volume = {54},
	year = {1932}}

@article{Hooker96,
	author = {Hooker, John N.},
	date = {1996/07/01},
	doi = {10.1007/BF02592142},
	id = {Hooker1996},
	isbn = {1436-4646},
	journal = {Mathematical Programming},
	number = {1},
	pages = {1--10},
	title = {Resolution and the integrality of satisfiability problems},
	url = {https://doi.org/10.1007/BF02592142},
	volume = {74},
	year = {1996}}

@article{Hooker88,
	author = {Hooker, John N.},
	date = {1988/03/01/},
	doi = {https://doi.org/10.1016/0167-9236(88)90097-8},
	isbn = {0167-9236},
	journal = {Decision Support Systems},
	number = {1},
	pages = {45--69},
	title = {A quantitative approach to logical inference},
	url = {https://www.sciencedirect.com/science/article/pii/0167923688900978},
	volume = {4},
	year = {1988}}

@article{Abdi-res,
	author = {Abdi, Ahmad and Cornu\'{e}jols, G\'{e}rard and Lee, Dabeen},
	journal = {Mathematics of Operations Research},
	number = {1},
	pages = {82-114},
	title = {Resistant Sets in the Unit Hypercube},
	volume = {46},
	year = {2020}}

@article{Fulkerson64,
	author = {D. Ray Fulkerson and O. A. Gross},
	date = {1964/9/1},
	journal = {Bulletin of the American Mathematical Society},
	journal1 = {Bulletin of the American Mathematical Society},
	journal2 = {Bulletin of the American Mathematical Society},
	month = {9},
	number = {5},
	pages = {681--684},
	title = {Incidence matrices with the consecutive 1's property},
	volume = {70},
	year = {1964}}

@incollection{Woodall78,
	author = {Woodall, Douglas R.},
	booktitle = {Theory and Applications of Graphs.},
	editor = {Alavi, Y. and Lick, D.R.},
	publisher = {Springer, Berlin, Heidelberg},
	series = {Lecture Notes in Mathematics},
	title = {Menger and {K}{\"o}nig systems},
	volume = {642},
	year = {1978}}

@article{Abdi-dijoins,
	author = {Ahmad Abdi and G{\'e}rard Cornu{\'e}jols and Michael Zlatin},
	journal = {SIAM Journal on Discrete Mathematics},
	number = {4},
	pages = {2417-2461},
	title = {On packing dijoins in digraphs and weighted digraphs},
	volume = {37},
	year = {2023}}

@article{Lucchesi78,
	author = {Lucchesi, Cl{\'a}udio L. and Younger, Daniel H.},
	doi = {10.1112/jlms/s2-17.3.369},
	fjournal = {Journal of the London Mathematical Society. Second Series},
	issn = {0024-6107},
	journal = {J. London Math. Soc. (2)},
	mrclass = {05C20},
	mrnumber = {500618},
	mrreviewer = {Jean H. Bevis},
	number = {3},
	pages = {369--374},
	title = {A minimax theorem for directed graphs},
	url = {https://doi-org.proxy.library.cmu.edu/10.1112/jlms/s2-17.3.369},
	volume = {17},
	year = {1978}}

@book{Lovasz86,
	author = {L{\'a}szl{\'o} Lov{\'a}sz and Michael D. Plummer},
	publisher = {North-Holland, Amsterdam},
	series = {North-Holland Mathematics Studies},
	title = {Matching Theory},
	volume = {121},
	year = {1986}}

@article{Shelah72,
	author = {Saharon Shelah},
	journal = {Pacific Journal of Mathematics},
	pages = {247--261},
	title = {A combinatorial problem; stability and order for models and theories in infinitary languages},
	volume = {41},
	year = {1972}}

@article{Sauer72,
	author = {Norbert Sauer},
	journal = {J. Comb. Theory, Ser. A},
	pages = {145--147},
	title = {On the density of families of sets},
	volume = {13},
	year = {1972}}

@article{Vapnik71,
	author = {Vladimir N. Vapnik and Alexey Y. Chervonenkis},
	journal = {Theory of Probability and its Applications},
	number = {2},
	pages = {264--280},
	title = {On the uniform convergence of relative frequencies of events to their probabilities},
	volume = {16},
	year = {1971}}

@article{Galvin14,
	author = {David Galvin},
	journal = {arXiv, 1406.7872},
	rss-description = {@misc{Galvin14,
      title={}, 
      author={David Galvin},
      year={2014},
      eprint={},
      archivePrefix={arXiv},
      primaryClass={math.CO},
      url={}, 
}},
	title = {Three tutorial lectures on entropy and counting},
	url = {https://arxiv.org/abs/1406.7872},
	year = {2014}}

@book{Diestel25,
	author = {Reinhard Diestel},
	publisher = {Springer Berlin, Heidelberg},
	title = {Graph Theory, 6th Edition},
	year = {2025}}

@article{Esperet11,
	author = {Esperet, Louis and Kardo{\v s}, Franti{\v s}ek and King, Andrew D. and Kr{\'a}l', Daniel and Norine, Serguei},
	journal = {Advances in Mathematics},
	number = {4},
	pages = {1646--1664},
	title = {Exponentially many perfect matchings in cubic graphs},
	volume = {227},
	year = {2011}}

@article{Barvinok13,
	author = {Barvinok, Alexander},
	date = {2013/02/01},
	doi = {10.1007/s00493-013-2870-9},
	id = {Barvinok2013},
	isbn = {1439-6912},
	journal = {Combinatorica},
	number = {1},
	pages = {1--10},
	title = {A bound for the number of vertices of a polytope with applications},
	url = {https://doi.org/10.1007/s00493-013-2870-9},
	volume = {33},
	year = {2013}}

@article{Guenin16,
	author = {Guenin, Bertrand},
	date = {2016//},
	doi = {10.1016/J.DAM.2015.10.035},
	id = {DBLP:journals/dam/Guenin16},
	journal = {Discret. Appl. Math.},
	pages = {122--132},
	title = {A survey on flows in graphs and matroids.},
	url = {https://doi.org/10.1016/j.dam.2015.10.035},
	volume = {209},
	year = {2016}}

@article{Seymour81,
	author = {Paul D. Seymour},
	journal = {Europ. J. Combinatorics},
	pages = {257-290},
	title = {Matroids and multicommodity flows},
	volume = {2},
	year = {1981}}

@article{Guenin01,
	author = {Guenin, Bertrand},
	fjournal = {Journal of Combinatorial Theory. Series B},
	issn = {0095-8956},
	journal = {J. Combin. Theory Ser. B},
	mrclass = {05C75},
	mrnumber = {1855799},
	mrreviewer = {Daniel Turz\'{\i}k},
	number = {1},
	pages = {112--168},
	title = {A characterization of weakly bipartite graphs},
	volume = {83},
	year = {2001}}

@book{Mitzenmacher17,
	address = {Cambridge},
	author = {Mitzenmacher, Michael and Upfal, Eli},
	db = {Cambridge Core},
	doi = {DOI: 10.1017/CBO9780511813603},
	dp = {Cambridge University Press},
	edition = {Second},
	publisher = {Cambridge University Press},
	title = {Probability and Computing: Randomized Algorithms and Probabilistic Analysis},
	url = {https://www.cambridge.org/core/product/3A5B47DB315FC64B9256C5C8131C5EFA},
	year = {2017}}

@incollection{Edmonds77,
	author = {Jack Edmonds and Rick Giles},
	booktitle = {Studies in Integer Programming},
	doi = {10.1016/S0167-5060(08)70734-9},
	editor = {P.L. Hammer and E.L. Johnson and B.H. Korte and G.L. Nemhauser},
	pages = {185-204},
	publisher = {Elsevier},
	series = {Annals of Discrete Mathematics},
	title = {A Min-Max Relation for Submodular Functions on Graphs},
	volume = {1},
	year = {1977}}

@article{Barahona86,
	author = {Barahona, Francisco and Gr{\"o}tschel, Martin},
	journal = {J. Combin. Theory Ser. B},
	number = {1},
	pages = {40-62},
	title = {On the cycle polytope of a binary matroid},
	volume = {40},
	year = {1986}}

@article{Cornuejols16,
	author = {Cornu{\'e}jols, G{\'e}rard and Lee, Dabeen},
	date = {2018/11/01},
	doi = {10.1007/s10107-017-1226-4},
	id = {Cornu{\'e}jols2018},
	isbn = {1436-4646},
	journal = {Mathematical Programming},
	number = {1},
	pages = {467--503},
	title = {On some polytopes contained in the 0, 1 hypercube that have a small Chv{\'a}tal rank},
	url = {https://doi.org/10.1007/s10107-017-1226-4},
	volume = {172},
	year = {2018}}

@article{Abdi-kwise,
	author = {Abdi, Ahmad and Cornu{\'e}jols, G{\'e}rard and Huynh, Tony and Lee, Dabeen},
	journal = {Math. Programming, Series B},
	pages = {29-50},
	title = {Idealness of $k$-wise intersecting families},
	volume = {192},
	year = {2021}}

@article{Benchetrit18,
	author = {Benchetrit, Yohann and Fiorini, Samuel and Huynh, Tony and Weltge, Stefan},
	journal = {Math. Oper. Res.},
	number = {3},
	pages = {718-725},
	title = {Characterizing polytopes in the 0/1-cube with bounded {C}hv{\'a}tal-{G}omory rank},
	volume = {43},
	year = {2018}}

@article{Abdi-cuboids,
	author = {Ahmad Abdi and G{\'e}rard Cornu{\'e}jols and Nat{\'a}lia Guri{\v c}anov{\'a} and Dabeen Lee},
	doi = {https://doi.org/10.1016/j.jctb.2019.10.002},
	issn = {0095-8956},
	journal = {Journal of Combinatorial Theory, Series B},
	pages = {144 - 209},
	title = {Cuboids, a class of clutters},
	url = {http://www.sciencedirect.com/science/article/pii/S0095895619301054},
	volume = {142},
	year = {2020}}

@article{Edmonds73,
	author = {Edmonds, Jack and Johnson, Ellis L.},
	doi = {10.1007/bf01580113},
	issn = {1436-4646},
	journal = {Mathematical Programming},
	month = {Dec},
	number = {1},
	pages = {88--124},
	publisher = {Springer Science and Business Media LLC},
	title = {Matching, {E}uler tours and the {C}hinese postman},
	url = {http://dx.doi.org/10.1007/BF01580113},
	volume = {5},
	year = {1973}}

@article{Fulkerson71,
	author = {Fulkerson, D. Ray},
	journal = {Math. Program.},
	pages = {168-194},
	title = {Blocking and anti-blocking pairs of polyhedra},
	volume = {1},
	year = {1971}}

@article{Seymour79,
	author = {Seymour, Paul D.},
	doi = {10.1112/plms/s3-38.3.423},
	eprint = {https://academic.oup.com/plms/article-pdf/s3-38/3/423/4351892/s3-38-3-423.pdf},
	issn = {0024-6115},
	journal = {Proceedings of the London Mathematical Society},
	month = {05},
	number = {3},
	pages = {423-460},
	title = {On multi-colourings of cubic graphs, and conjectures of Fulkerson and Tutte},
	url = {https://doi.org/10.1112/plms/s3-38.3.423},
	volume = {38},
	year = {1979}}

@book{Cornuejols01,
	author = {Cornu{\'e}jols, G{\'e}rard},
	publisher = {SIAM},
	title = {Combinatorial Optimization: Packing and Covering},
	volume = {74},
	year = {2001}}

@article{Cornuejols94,
	author = {Cornu\'{e}jols, G\'{e}rard and Novick, Beth},
	fjournal = {Journal of Combinatorial Theory. Series B},
	issn = {0095-8956},
	journal = {J. Combin. Theory Ser. B},
	mrclass = {05B20},
	mrnumber = {1256587},
	mrreviewer = {Ding Zhu Du},
	number = {1},
	pages = {145--157},
	title = {Ideal {$0,1$} matrices},
	volume = {60},
	year = {1994}}

@article{Edmonds70,
	author = {Edmonds, Jack and Fulkerson, D. R.},
	fjournal = {Journal of Combinatorial Theory},
	issn = {0021-9800},
	journal = {J. Combinatorial Theory},
	mrclass = {90.50 (05.00)},
	mrnumber = {255235},
	mrreviewer = {M. Dragomirescu},
	pages = {299--306},
	title = {Bottleneck extrema},
	volume = {8},
	year = {1970}}

@article{Isbell58,
	author = {Isbell, John R.},
	fjournal = {Duke Mathematical Journal},
	issn = {0012-7094},
	journal = {Duke Math. J.},
	mrclass = {52.00},
	mrnumber = {97028},
	mrreviewer = {J. H. Blau},
	pages = {423--439},
	title = {A class of simple games},
	volume = {25},
	year = {1958}}

@article{Lehman79,
	author = {Lehman, Alfred},
	fjournal = {Mathematical Programming},
	issn = {0025-5610},
	journal = {Math. Programming},
	mrclass = {94C15 (05C50 90B10)},
	mrnumber = {527577},
	mrreviewer = {Fran\c{c}ois Aribaud},
	number = {2},
	pages = {245--259},
	title = {On the width-length inequality},
	volume = {16},
	year = {1979}}
 
\end{document}